\documentclass[sigconf]{acmart}
\AtBeginDocument{%
  \providecommand\BibTeX{{%
    \normalfont B\kern-0.5em{\scshape i\kern-0.25em b}\kern-0.8em\TeX}}}

\copyrightyear{2024}
\acmYear{2024}
\setcopyright{rightsretained}
\acmConference[SPAA '24]{Proceedings of the 36th ACM Symposium on
Parallelism in Algorithms and Architectures}{June 17--21, 2024}{Nantes, France}
\acmBooktitle{Proceedings of the 36th ACM Symposium on Parallelism in
Algorithms and Architectures (SPAA '24), June 17--21, 2024, Nantes,
France}\acmDOI{10.1145/3626183.3659980}
\acmISBN{979-8-4007-0416-1/24/06}


\usepackage{fancyhdr}

\usepackage{algorithm}
\usepackage[noend]{algorithmic}

\usepackage{amsmath, amsthm}
\usepackage{graphicx}
\usepackage{accents}
\usepackage{hyperref}
\usepackage{bbm}
\usepackage{bm}
\usepackage{booktabs}
\usepackage{float}
\usepackage{multirow}
\usepackage{scalerel}
\newtheorem{theorem}{Theorem}[section]

\usepackage{tikz}
\usetikzlibrary{calc, decorations.pathmorphing}
\usetikzlibrary{decorations.pathreplacing, chains}
\usetikzlibrary {arrows.meta}

\newcommand{\RR}{\mathbb{R}}

\newcommand{\ZZ}{\mathbb{Z}}

\newcommand{\paren}[1]{\left( #1 \right)}
\newcommand{\br}[1]{\left[ #1 \right]}

\newcommand{\norm}[1]{\left\lVert#1\right\rVert}

\newcommand{\scr}[1]{\mathcal{#1}}
\newcommand{\nnz}[1]{\textrm{nnz}(#1)}

\newcommand{\pluseq}{\mathrel{+}=}

\newcommand{\diveq}{\mathrel{/}=}
\newcommand{\argmin}{\textrm{argmin}}
\mathchardef\mhyphen="2D

\DeclareMathOperator*{\startimes}{\scalerel*{\circledast}{\sum}}

\begin{document}

\title{Distributed-Memory Randomized Algorithms for Sparse 
Tensor CP Decomposition}

\author{Vivek Bharadwaj}
\affiliation{%
  \institution{University of California, Berkeley}
  \city{Berkeley}
  \state{CA}
  \country{USA}
}
\email{vivek_bharadwaj@berkeley.edu}

\author{Osman Asif Malik}
\authornote{Work completed while author was affiliated with
Lawrence Berkeley National Laboratory.}
\affiliation{%
  \institution{Encube Technologies}
  \city{Stockholm}
  \country{Sweden}
}
\email{osman@getencube.com}

\author{Riley Murray}
\authornote{Work completed while author was affiliated
with Lawrence Berkeley National Laboratory,
International Computer Science Institute, 
and University of California, Berkeley.}
\affiliation{%
  \institution{Sandia National Laboratories}
  \city{Albuquerque}
  \state{NM}
  \country{USA}
}

\email{rjmurr@sandia.gov}

\author{Aydın Buluç}
\affiliation{%
  \institution{Lawrence Berkeley National Laboratory}
  \city{Berkeley}
  \state{CA}
  \country{USA}
}
\email{abuluc@lbl.gov}

\author{James Demmel}
\affiliation{%
  \institution{University of California, Berkeley} 
  \city{Berkeley}
  \state{CA}
  \country{USA}
}
\email{demmel@berkeley.edu}

\renewcommand{\shortauthors}{Vivek Bharadwaj,
Osman Asif Malik, Riley Murray, Aydın Buluç, \& James Demmel
}

\begin{abstract}
Candecomp / PARAFAC (CP) decomposition, a generalization of 
the matrix singular value 
decomposition to higher-dimensional tensors, is a popular tool for analyzing 
multidimensional sparse data. On tensors with billions of nonzero entries, 
computing a CP decomposition is a computationally intensive task. 
We propose the first distributed-memory implementations
of two randomized CP decomposition algorithms,
CP-ARLS-LEV and STS-CP, that offer nearly an order-of-magnitude
speedup at high decomposition ranks over well-tuned non-randomized decomposition 
packages. Both algorithms rely on leverage score sampling and enjoy 
strong theoretical guarantees, each with varying time and 
accuracy tradeoffs. 
We tailor the communication schedule for our 
random sampling algorithms, 
eliminating expensive reduction collectives and forcing communication costs to scale with the random sample 
count. 
Finally, we optimize the local storage format for our methods, switching between analogues of compressed sparse column and compressed sparse row 
formats. Experiments show that our methods are fast and scalable,
producing 11x speedup over SPLATT by decomposing the billion-scale 
Reddit tensor on 512 CPU cores in under two minutes. 
\end{abstract}



\maketitle

\section{Introduction}
\begin{figure}[h]
    \centering
    \includegraphics[scale=0.2]{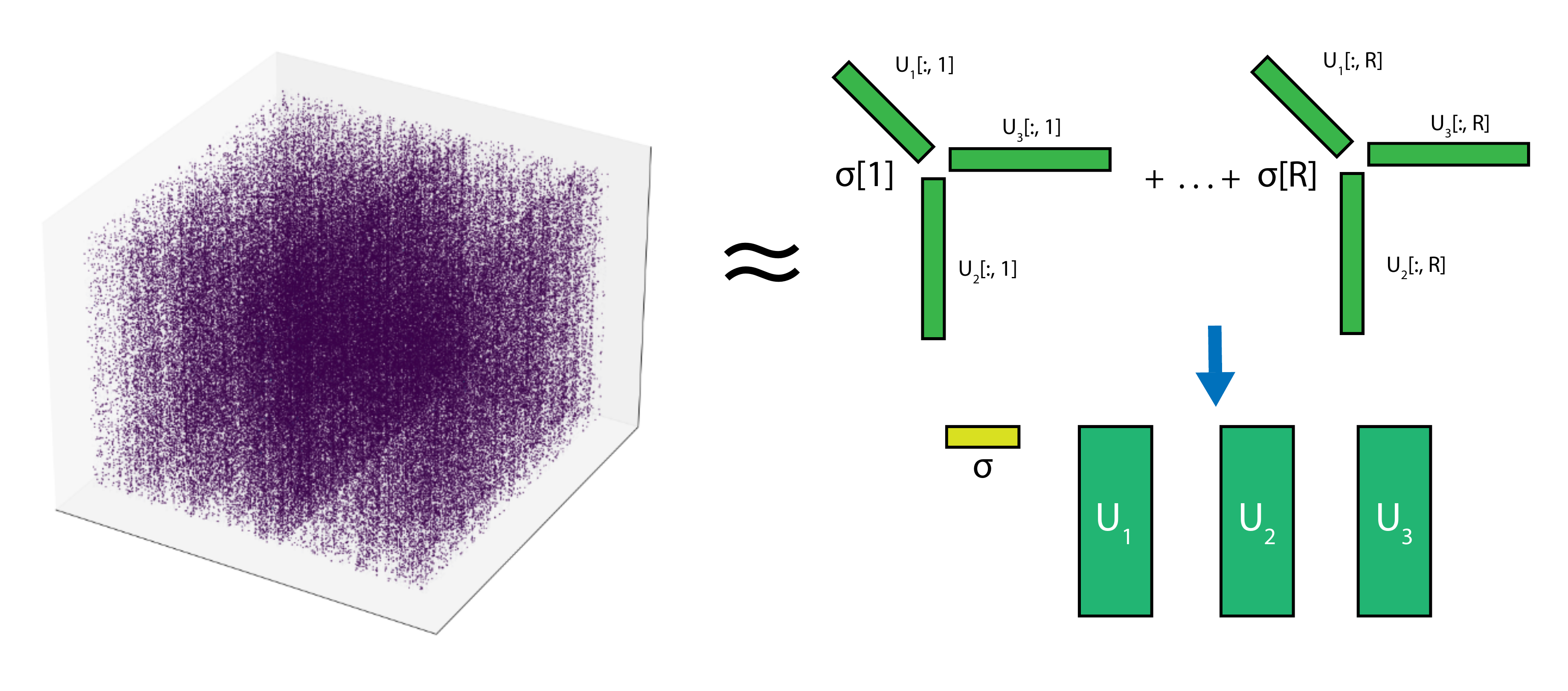}
    \caption{A subset of entries from the 3D
    Amazon Review sparse tensor \cite{smith_frostt_2017} 
    and its illustrated CP decomposition. The
    sparse tensor is approximated by a sum
    of 3D outer products of vectors, which are 
    columns of factor matrices $U_1, U_2$, and $U_3$. 
    Outer products are scaled by elements of $\sigma$.}
    \Description[CP Decomposition Illustration]{Left: An illustration of a cube with several
    discrete points representing the sparse tensor. Right: outer products,
    summed together, each of which is represented
    by a triplet of vectors radiating from
    a central point.} 
    \label{fig:cp-illustration}
\end{figure}

Randomized algorithms for numerical
linear algebra have become increasingly popular 
in the past decade, but their distributed-memory
communication characteristics and scaling properties 
have received less attention. In this work, we examine 
randomized algorithms to compute
the Candecomp / PARAFAC (CP) decomposition, a
generalization of the matrix singular-value
decomposition to a number of modes $N > 2$.
Given a tensor $\scr T \in  
\RR^{I_1 \times ... \times I_N}$ and a
target rank $R$, the goal of CP decomposition (illustrated in Figure
\ref{fig:cp-illustration}) is to find
a set of \textit{factor matrices} $U_1, ...,
U_N, U_j \in \RR^{I_j \times R}$ with unit norm columns
and a nonnegative vector $\sigma \in \RR^R$ satisfying 
\begin{equation}    
\scr T\br{i_1, ..., i_N} 
\approx \sum_{r=1}^R \sigma\br{r} * 
U_1\left[i_1, r\right] * ... * U_N\br{i_N, r}.
\label{eq:cp_decomposition_defn}
\end{equation}
We consider real \textbf{sparse} tensors $\scr T$ 
with $N \geq 3$, all entries known, and billions of nonzero entries. 
Sparse tensors are a flexible abstraction
for a variety of data, such as network traffic logs \cite{mao_malspot_2014}, 
text corpora \cite{smith_frostt_2017}, and knowledge graphs \cite{balazevic_tucker_2019}.

\subsection{Motivation}
Why is a low-rank 
approximation of 
a sparse tensor useful? We can view the 
sparse CP decomposition as an extension of 
well-studied sparse matrix factorization
methods, which can 
mine patterns from large datasets \cite{haesun_matrix_nmf}.
Each row of the CP factors is a dense embedding 
vector for an index 
$i_j \in \br{I_j}$, $1 \leq j \leq N$.
Because each embedding is a small
dense vector while the input 
tensor is sparse, 
sparse tensor CP decomposition may 
incur high relative error with respect
to the input and rarely captures the
tensor sparsity structure exactly.
Nevertheless, the learned embeddings 
contain valuable information. 
CP factor matrices have 
been successfully used to identify 
patterns in social 
networks \cite{gcp_kolda, larsen_practical_2022},
detect anomalies in packet 
traces \cite{mao_malspot_2014}, and monitor
trends in internal network traffic 
\cite{smith_streaming}. As we discuss below,
a wealth of software packages exist to meet the demand
for sparse tensor decomposition.

One of the most popular methods for computing
a sparse CP decomposition, the 
\textit{Alternating-Least-Squares} (ALS) algorithm, 
involves repeatedly solving large, 
overdetermined linear least-squares problems
with structured design matrices \cite{kolda_tensor_2009}. 
High-performance libraries DFacto\cite{choi_dfacto_2014},
SPLATT \cite{smith_splatt_2015},
HyperTensor \cite{hypertensor}, 
and BigTensor \cite{bigtensor_16} distribute these expensive 
computations to a cluster 
of processors that communicate through an interconnect. 
Separately, several works use randomized sampling methods to accelerate the
least-squares solves, 
with prototypes implemented in a shared-memory setting
\cite{cheng_spals_2016, larsen_practical_2022, malik_more_2022, bharadwaj2023fast}.
These randomized algorithms have strong theoretical guarantees and offer significant 
asymptotic advantages over non-randomized ALS. Unfortunately, prototypes of
these methods require hours to run \cite{larsen_practical_2022, bharadwaj2023fast}
and are neither competitive nor scalable compared to 
existing libraries with distributed-memory parallelism.

\begin{figure}
    \centering
    \includegraphics[scale=0.34]{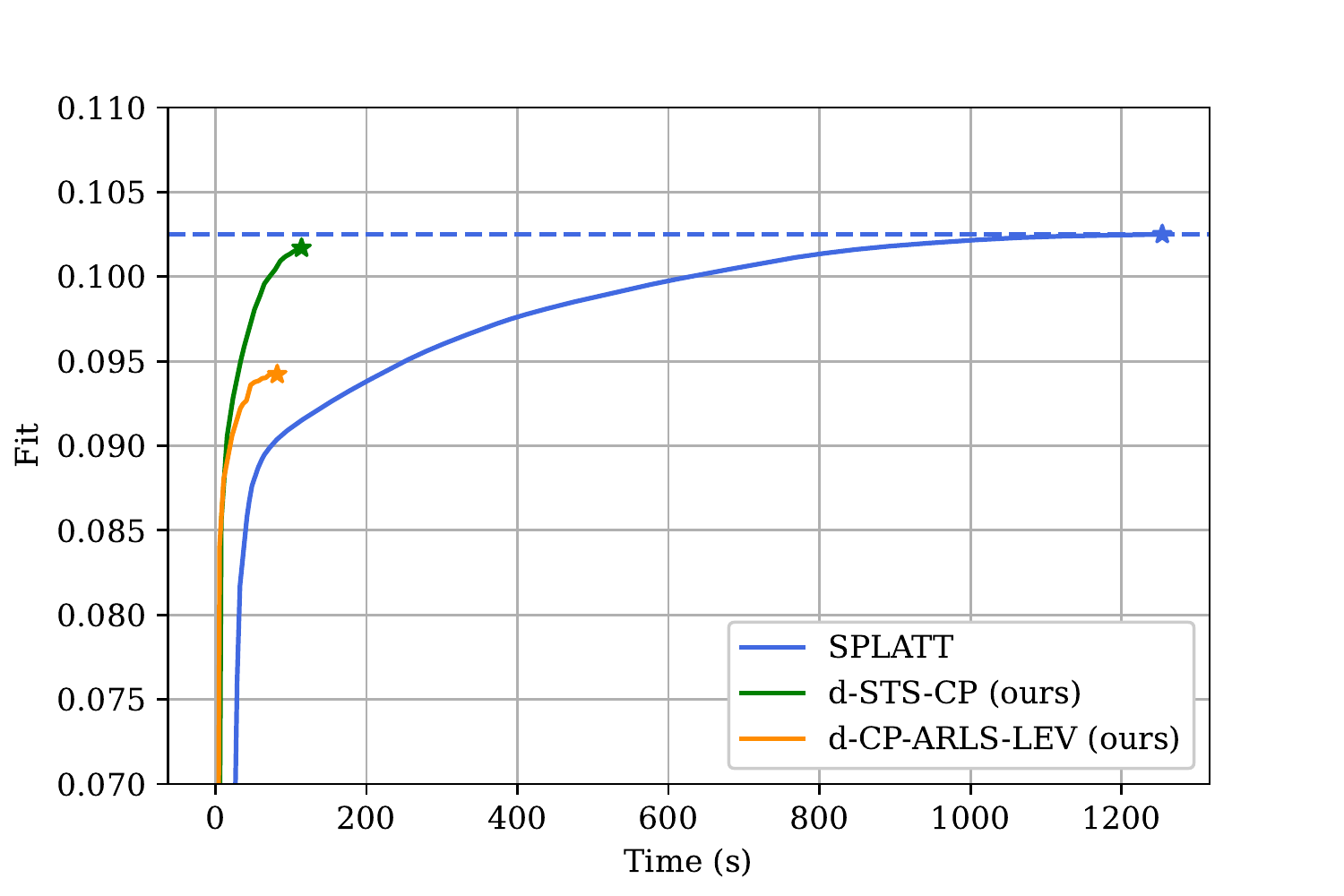}
    \caption{Running maximum accuracy over time for SPLATT, a state-of-the-art distributed
    CP decomposition software package, and
    our randomized algorithms on the Reddit 
    tensor, target rank $R=100$, on  
    512 CPU cores. Curves are averages 
    of 5 trials, 80 ALS rounds.}
    \Description[STS-CP and CP-ARLS-LEV make 
    faster progress
    to the solution]{Three curves rising from
    the bottom left corner to the top right. The
    first two correspond to CP-ARLS-LEV and 
    STS-CP, approaching maximum 
    accuracy within 120
    seconds. The third curve (SPLATT) for
    the non-randomized algorithm requires more 
    than 1200 seconds to approach its maximum. 
    The maximums for STS-CP and SPLATT are
    nearly identical.}
    \label{fig:reddit_fit_vs_time}
\end{figure}

\subsection{Our Contributions}
We propose the first distributed-memory parallel formulations 
of two randomized algorithms,
CP-ARLS-LEV \cite{larsen_practical_2022} and STS-CP
\cite{bharadwaj2023fast}, with accuracy identical to
their shared-memory prototypes.
We then provide implementations of these methods 
that scale to thousands of CPU cores. We face 
\textbf{dual technical challenges} to
parallel scaling. First, sparse tensor decomposition generally has 
lower arithmetic intensity (FLOPs / data word communicated between processors) than 
dense tensor decomposition, since computation scales linearly with the 
tensor nonzero count. Some sparse tensors exhibit nonzero fractions as low as 
$4 \times 10^{-10}$ (see Table \ref{tab:datasets}), while
the worst-case communication costs for sparse CP decomposition 
remain identical to the dense tensor case 
\cite{smith_medium-grained_2016}. Second, randomized
algorithms can save an order of magnitude in computation over their 
non-randomized counterparts 
\cite{mahoney_rand_survey, drineas_mahoney_rnla_survey, martinsson_tropp_2020}, but their inter-processor communication costs
remain unaltered unless carefully optimized. Despite these compounding factors
that reduce arithmetic intensity, we achieve both speedup and 
scaling through several key innovations, three of which we highlight:

    \paragraph{Novel Distributed-Memory
    Sampling Procedures} Random sample selection 
    is challenging to implement when the CP factor matrices
    and sparse tensor are divided among $P$ processors.
    We introduce two distinct communication-avoiding algorithms for 
    randomized sample selection from the 
    Khatri-Rao product. First, we show how to implement
    the CP-ARLS-LEV algorithm by computing an 
    independent probability
    distribution on the factor block row owned 
    by each processor. The resulting distributed
    algorithm has minimal compute / communication 
    overhead compared to the other phases of CP
    decomposition. The second algorithm, STS-CP, requires
    higher sampling time, but achieves lower error by
    performing random walks on a binary tree for each sample.
    By distributing leaf nodes uniquely to processors and replicating
    internal nodes, we give a sampling algorithm 
    with per-processor communication bandwidth scaling as 
    $O\paren{\log P / P}$ (see Table \ref{tab:sampler_compute_comm}). 

    \paragraph{Communication-Optimized MTTKRP} 
    We show that 
    communication-optimal schedules for non-randomized ALS may exhibit disproportionately 
    high communication costs for randomized algorithms.
    To combat this, we use an ``accumulator-stationary" schedule that eliminates
    expensive \verb|Reduce-scatter| collectives, causing all communication costs
    to scale with the number of random samples taken. This alternate schedule significantly reduces 
    communication on tensors with large 
    dimensions (Figure 
    \ref{fig:communication_comparison}) 
    and empirically improves the computational
    load balance (Figure \ref{fig:load_imbalance}).

    \paragraph{Local Tensor Storage Format} Existing
    storage formats developed for sparse
    CP decomposition \cite{smith_splatt_2015, 
    nisa_csf_mm} are not optimized for random access into
    the sparse tensor, which our algorithms require. 
    In response, we use a modified compressed-sparse-column format to store each matricization of our tensor, allowing efficient selection
    of nonzero entries by our random sampling algorithms. We then transform the
    selected nonzero entries into compressed sparse row format,
    which eliminates shared-memory data races in the subsequent sparse-dense matrix
    multiplication. The cost of the transposition is justified and provides a roughly 1.7x speedup over using atomics in a
    hybrid OpenMP / MPI implementation.
    \newline

Our distributed-memory randomized algorithms have significant advantages over existing libraries while preserving the accuracy
of the final approximation. As Figure 
\ref{fig:reddit_fit_vs_time} shows, our method d-STS-CP computes a rank 100 
decomposition of the Reddit tensor ($\sim 4.7$ billion nonzero entries) with a 11x
speedup over SPLATT, a state-of-the-art distributed-memory CP decomposition 
package. The reported speedup was achieved on 512 CPU cores, with a final fit within $0.8\%$ of non-randomized ALS for the same iteration count. 
While the distributed algorithm d-CP-ARLS-LEV 
achieves a lower final accuracy, it makes progress faster than
SPLATT and spends less time on sampling (completing 80 
rounds in an average of 81 seconds). We demonstrate that 
it is well-suited to smaller tensors and lower target ranks.

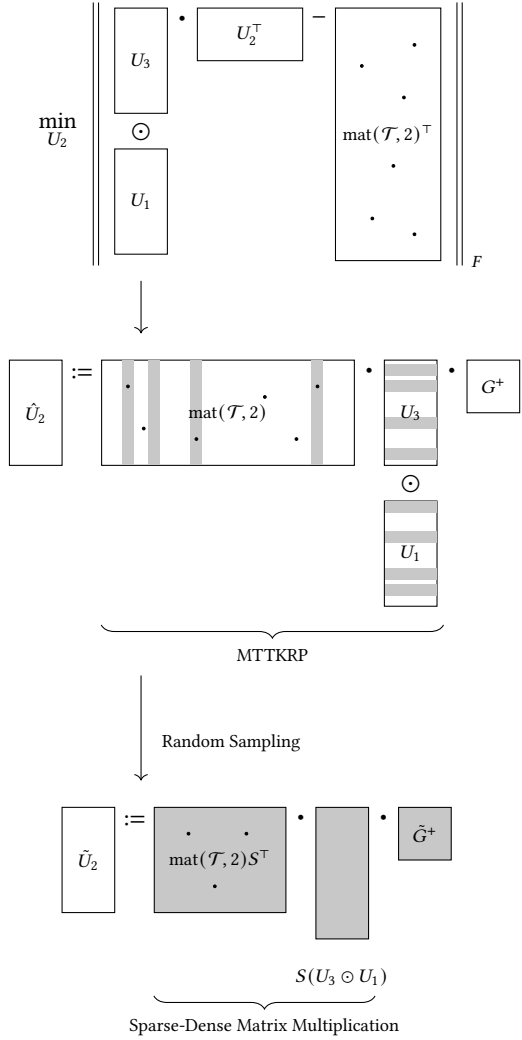
\begin{figure}
    \centering
    \begin{tikzpicture}[scale=0.7,
mnd/.style n args={2}{transform shape, draw, 
minimum height=#1 cm, minimum width=#2 cm,anchor=center,font=\footnotesize}
]
%\draw[help lines] grid (9,9);
\node[mnd={2}{1}] at (2, 8) (U3) {};
\node[] at (U3.center) {\footnotesize $U_3$};
\node[below, font=\Large] at (U3.south) (odot) {$\odot$};
\node[right, anchor=north west, font=\Huge] at (U3.north east) (cdot) {$\cdot$};
\node[below, mnd={2}{1}, anchor=north] at (odot.south) (U1) {};
\node[] at (U1.center) {\footnotesize $U_1$};
\node[right, mnd={1}{2}, anchor=north west] at (cdot.north east) (U2) {};
\node[] at (U2.center) {\footnotesize $U_2^\top$};
\node[right, anchor=north west, font=\normalsize] at (U2.north east) (minus) {$-$};
\node[right, mnd={4.8}{2}, anchor=north west] at (minus.north east) (spmat) {};
\newcommand{\raggedheight}{2cm}
\foreach \val in {(-0.5, 1.3), (-0.3, -1.6), (0.1, -0.6), (0.3, 0.7), (0.5, 1.7), (0.5, -1.9)}
{
\draw[fill=black] (spmat.center) ++\val circle(0.08em);
}
\node[font=\footnotesize] at (spmat.center) {$\textrm{mat}(\scr T, 2)^\top$};
\draw (U3.north west) ++ (-0.3, 0.1) -- ++ (0, -5cm) ;
\draw (U3.north west) ++ (-0.4, 0.1) -- ++ (0, -5cm) ;
\draw (spmat.north east) ++ (0.3, 0.1) -- ++ (0, -5cm) ;
\draw (spmat.north east) ++ (0.4, 0.1) -- ++ (0, -5cm) ;
\node[xshift=-3.5em] at (odot.center) {$\min\limits_{U_2}$};
\node[xshift=1.5em, font=\scriptsize] at (spmat.south east) {$F$};

% Second image
\draw[->] (U1.south) ++(0, -0.5)  -- ++(0, -1); 
\node[mnd={2}{1}, anchor=north, xshift=-2cm, yshift=-2.0cm] at (U1.south) (U2) {};
\node[] at (U2.center) {\footnotesize $\hat U_2$};
\node[right, anchor=north west, font=\large] at (U2.north east) (ceq) {$:=$};
\node[right, mnd={2}{4.8}, anchor=north west] at (ceq.north east) (spmat) {};

\definecolor{gray_cust}{RGB}{200,200,200} 
% Gray shading showing rows taken from the MTTKRP

\foreach \x in {1.0cm, 0.5cm, 1.8cm, 4.1cm}
{
\node[mnd={1.98}{0.2}, anchor=center, fill=gray_cust, draw=none,
xshift=\x] at (spmat.west) {};
}

% Put points in the sparse tensor
\foreach \val in {(1.3, -0.5), (-1.6, -0.3), (-0.6, -0.5), 
(0.7, 0.3), (1.7, 0.5), (-1.9, 0.5)}
{
\draw[fill=black] (spmat.center) ++\val circle(0.08em);
}
\node[font=\footnotesize] at (spmat.center) {$\textrm{mat}(\scr T, 2)$};
\node[right, anchor=north west, font=\Huge] at (spmat.north east) (cdot1) {$\cdot$};
\node[mnd={2}{1}, anchor=north west] at (cdot1.north east) (U3) {};

\foreach \y in {-0.2cm, -0.5cm, -1.2cm, -1.8cm}
{
\node[mnd={0.2}{0.98}, anchor=center, fill=gray_cust, draw=none,
yshift=\y] at (U3.north) {};
}
\node[] at (U3.center) {\footnotesize $U_3$};

\node[below, font=\Large] at (U3.south) (odot) {$\odot$};
\node[right, anchor=north west, font=\Huge] at (U3.north east) (cdot2) {$\cdot$};
\node[below, mnd={2}{1}, anchor=north] at (odot.south) (U1_new) {};

\foreach \y in {-0.13cm, -0.7cm, -1.4cm, -1.7cm}
{
\node[mnd={0.2}{0.98}, anchor=center, fill=gray_cust, draw=none,
yshift=\y] at (U1_new.north) {};
}

\node[] at (U1_new.center) {\footnotesize $U_1$};
\node[below, mnd={1}{1}, anchor=north west] at (cdot2.north east) (Gpinv) {};
\node[] at (Gpinv.center) {\footnotesize $G^+$};
\draw [decorate,decoration={brace,amplitude=5pt,mirror,raise=4ex}]
  (spmat.south west) ++ (0, -2.25)  -- ++(6.5,0) node[midway,yshift=-3em]{\footnotesize MTTKRP};

% Third image
\draw[->] (U1.south) ++(0, -8.0)  -- ++(0, -2); 
\node[xshift=1.2cm, yshift=-6.5cm] at (U1.south) {\footnotesize Random Sampling};

\node[mnd={2}{1}, anchor=north, xshift=-1cm, yshift=-10.5cm] at (U1.south) (U2) {};
\node[] at (U2.center) {\footnotesize $\tilde U_2$};
\node[right, anchor=north west, font=\large] at (U2.north east) (ceq) {$:=$};
\node[right, mnd={2}{2.5}, anchor=north west, fill=gray_cust] at (ceq.north east) (spmat) {};

% Put points in the sparse tensor
\foreach \val in {(-0.6, 0.5), (0.5, 0.5), (-0.1, -0.5)}
{
\draw[fill=black] (spmat.center) ++\val circle(0.08em);
}
\node[font=\footnotesize] at (spmat.center) {$\textrm{mat}(\scr T, 2) S^\top$};
\node[right, anchor=north west, font=\Huge] at (spmat.north east) (cdot1) {$\cdot$};
\node[mnd={2.5}{1}, anchor=north west, fill=gray_cust] at (cdot1.north east) (U3) {};
\node[yshift=-0.5cm] at (U3.south) {\footnotesize $S (U_3 \odot U_1)$};
\node[right, anchor=north west, font=\Huge] at (U3.north east) (cdot2) {$\cdot$};
\node[below, mnd={1}{1}, anchor=north west, fill=gray_cust] at (cdot2.north east) (Gpinv) {};
\node[] at (Gpinv.center) {\footnotesize $\tilde G^+$};
\draw [decorate,decoration={brace,amplitude=5pt,mirror,raise=4ex}]
  (spmat.south west) ++ (0, -0.80)  -- ++(4.2,0) node[midway,yshift=-3em]{\footnotesize Sparse-Dense Matrix Multiplication};

\end{tikzpicture}
    \caption{Top: the linear least-squares problem to optimize factor matrix $U_2$ during the ALS algorithm for a 3D tensor (column dimension of 
    $\textrm{mat}(\scr T, 2)$ not to scale). Middle: the exact 
    solution to the problem using the Matricized Tensor 
    Times Khatri-Rao Product
    (MTTKRP). Shaded columns of $\textrm{mat}(\scr T, 2)$ and rows
    of $(U_3 \odot U_1)$ are selected by our 
    random sampling algorithm. Bottom: the downsampled
    linear least-squares problem after applying random 
    sampling matrix $S$.}
    \label{fig:mttkrp_definition}
    \Description[The Randomized MTTKRP Solve
    Procedure]{Top: illustration of the main ALS
    equation solving for U2 with tall-skinny
    matrices U3 and U1, along with a tall 
    rectangle for the matricized tensor with
    dots inside indicating its sparsity. Middle:
    the normal equation solution, with columns
    from the matricized tensor and rows from U3, 
    U1 shaded in gray. Bottom: a normal
    equation illustration with much smaller
    matrices.}
\end{figure}

\section{Notation and Preliminaries}

\begin{table}
\centering
\begin{tabular}{ll}
\toprule
\textbf{Symbol} & \textbf{Description} \\
\midrule
$\scr T$ & Sparse tensor of dimensions $I_1 \times \ldots \times I_N$ \\
$R$ & Target Rank of CP Decomposition \\
%$\textrm{nnz}(\scr S)$ & Number of nonzeros in $\scr S$ \\
$U_1, \ldots, U_N$ & Dense factor matrices, $U_j \in \RR^{I_j \times R}$ \\
$\sigma$ & Vector of scaling 
factors, $\sigma \in \RR^{R}$\\
$J$ & Sample count for randomized ALS \\
\midrule
$\cdot$ & Matrix multiplication \\
$\circledast$ & Elementwise multiplication \\
$\otimes$ & Kronecker product \\
$\odot$ & Khatri-Rao product \\
\midrule
$P$ & Total processor count \\
$P_1, \ldots, P_N$ & Dimensions of processor grid, $\prod_i P_i = P$ \\
$U_i^{(p_j)}$ & Block row of $U_i$ owned by processor $p_j$ \\
\bottomrule
\end{tabular}
\vspace{5pt}
\caption{Symbol Definitions}
\label{table:symbol_defs}
\end{table}

Table \ref{table:symbol_defs} summarizes our notation. We use script 
characters (e.g. $\scr T$) to denote tensors with at 
least three modes, capital letters for matrices, and lowercase letters for vectors. Bracketed 
tuples following any of these objects, e.g. $A\br{i, j}$, 
represent indexes into each object, and the 
symbol ``:" in place of any index 
indicates a slicing of a tensor. We use 
$\odot$ to denote the Khatri-Rao product, which is a column-wise Kronecker product of a pair of matrices with the 
same number of columns. For $A \in \RR^{I \times R}, 
B \in \RR^{J \times R}$, $A \odot B$ produces a matrix of dimensions $(IJ) \times R$ such that for $1 \leq j \leq R$,  
$$(A \odot B)\br{:, j} = A\br{:, j} \otimes B\br{:, j}.$$

Let $\scr T$ be an $N$-dimensional tensor indexed by tuples 
$(i_1, ..., i_N) \in \br{I_1} \times ... \times \br{I_N}$, with $\nnz{\scr T}$ as the number of
nonzero entries. In this work, sparse tensors 
are always represented as a collection
of $(N+1)$-tuples, with the first $N$ elements giving the indices of a nonzero element
and the last element giving the value.
We seek a low-rank approximation of $\scr T$ given by Equation 
\eqref{eq:cp_decomposition_defn}, the right-hand-side of which
we abbreviate as $\br{\sigma; U_1, ..., U_N}$. 
By convention, each column of $U_1, ..., U_N$ has unit norm. Our goal is 
to minimize the sum of squared differences between our approximation 
and the provided tensor: 
\begin{equation}
\textrm{argmin}_{\sigma, U_1, \ldots, U_N} \norm{\br{\sigma; U_1, ..., U_N} - \scr T}_{\textrm{F}}^2.
\label{eq:frob_norm}
\end{equation}

\subsection{Non-Randomized ALS CP Decomposition}
Minimizing Equation \eqref{eq:frob_norm} jointly over $U_1, ..., U_N$ is still a non-convex problem (the vector $\sigma$ can be
computed directly from the factor matrices
by renormalizing each column). Alternating least squares is a popular heuristic algorithm that iteratively drives down the approximation error. The algorithm begins with a set 
of random factor matrices and optimizes the approximation in rounds, 
each involving $N$ subproblems. The $j$-th
subproblem in a round holds all factor matrices but $U_j$ constant and solves for a new matrix 
$\hat U_j$ minimizing the squared Frobenius norm error \cite{kolda_tensor_2009}. The updated matrix $\hat U_j$ is the solution to the overdetermined
linear least-squares problem
\begin{equation}
\hat U_j := \min_{X} \norm{U_{\neq j} \cdot X^\top 
- \textrm{mat}(\scr T, j)^\top}_{F}.
\label{eq:cp_decomposition}  
\end{equation}
Here, the design matrix is 
\[U_{\neq j} := U_N \odot ... \odot U_{j+1} \odot U_{j-1} \odot ... \odot U_1,
\]
which is a Khatri-Rao Product (KRP) of the factors held constant.
The matrix $\textrm{mat}(\scr T, j)$ is a 
\textit{matricization} of the sparse tensor $\scr{T}$, which
reorders the tensor modes and flattens it into a matrix of dimensions
$I_j \times (\prod_{i \neq j} I_i)$. We solve the problem 
efficiently using the normal equations. Denoting the Gram 
matrix by $G = (U_{\neq j})^\top (U_{\neq j})$, we have 
\begin{equation} 
    \hat U_j := \textrm{mat}(\scr T, j) \cdot U_{\neq j} \cdot G^+,
    \label{eq:cp_update_normal} 
\end{equation}
where $G^+$ is the Moore-Penrose pseudo-inverse of $G$.
Since $U_{\neq j}$ is a Khatri-Rao product, we 
can efficiently compute $G$ through the 
well-known \cite{kolda_tensor_2009} formula 
\begin{equation}
G = \startimes_{k \neq j} (U_k^\top U_k),
\label{eq:efficient_G_formula}
\end{equation}    
where $\circledast$ denotes elementwise multiplication.
Figure \ref{fig:mttkrp_definition} illustrates each
least-squares problem, and Algorithm \ref{alg:als_algorithm} 
summarizes the ALS procedure, including a renormalization of factor matrix columns after each solve. We implement the 
initialization step in line 1 by drawing all factor matrix entries from a 
unit-variance Gaussian distribution, a standard 
technique \cite{larsen_practical_2022}.

\begin{algorithm}
   \caption{CP-ALS($\scr T$, $R$)} 
   \label{alg:als_algorithm}
\begin{algorithmic}[1]
        \STATE Initialize $U_j \in \RR^{I_j \times R}$
        randomly for $1 \leq j \leq N$.

        \STATE Renorm. $U_j\br{:, i} \diveq \norm{U_j\br{:, i}}_2 $, 
        $1 \leq j \leq N, 1 \leq i \leq R$.

        \STATE Initialize $\sigma \in \RR^{R}$ to $\br{1}$. 

        \WHILE{not converged}
            \FOR{$j=1...N$}
                \STATE $U_j := \argmin_X \norm{U_{\neq j} \cdot X^\top 
                -  \textrm{mat}(\scr T, j)^\top}_F $

                \STATE $\sigma \br{i} = \norm{U_j \br{:, i}}_2$, $1 \leq i \leq R$

                \STATE Renorm. $U_j\br{:, i} \diveq \norm{U_j\br{:, i}}_2 $, 
                $1 \leq i \leq R$.
            \ENDFOR
        \ENDWHILE

        \STATE \textbf{return} $\br{\sigma; U_1, ..., U_N}$.
\end{algorithmic}
\end{algorithm}
The most expensive component of the ALS algorithm 
is the matrix multiplication 
$\textrm{mat}(\scr T, j) \cdot U_{\neq j}$ in Equation \eqref{eq:cp_update_normal}, 
an operation known as the Matricized Tensor-Times-Khatri Rao Product 
(MTTKRP). For a sparse tensor $\scr T$, this kernel has a 
computational pattern similar to sparse-dense 
matrix multiplication (SpMM): for each nonzero in the sparse 
tensor, we compute a scaled Hadamard
product between $N-1$ rows of the constant factor matrices and add it to a row of the remaining factor matrix. The MTTKRP runtime is 
\begin{equation}
 O(\nnz{\scr T} N R),
\label{eq:mttkrp_runtime}
\end{equation}
which is linear in the nonzero count of $\scr T$. Because $\scr T$ may 
have billions of nonzero entries, we seek methods to
drive down the cost of the MTTKRP. 

\subsection{Randomized Leverage Score Sampling}
Sketching is a powerful tool to accelerate least squares 
problems of the form $\min_X\|AX-B\|_{\text{F}}$ where $A$ has far more rows than columns 
\cite{mahoney_rand_survey, drineas_mahoney_rnla_survey, martinsson_tropp_2020}.
We apply a structured 
sketching matrix $S^{J \times I}$ to both
$A$ and $B$, where the
row count of $S$ satisfies $J \ll I$.
The resulting problem $\min_{\tilde X}\|S(A\tilde X-B)\|_{F}$ is cheaper to solve,
and the solution $\tilde X$ has residual arbitrarily close (for sufficiently high $J$) to the true minimum with high probability. We seek a sketching operator $S$ with an efficiently computable action on $A$, which is a
Khatri-Rao product.

We choose $S$ to be a \textit{sampling} matrix with a single nonzero per row (see
Section \ref{sec:alternate_sketches} 
for alternatives). This matrix 
extracts and reweights $J$ rows from both $A$ and
$B$, preserving the 
sparsity of the matricized tensor $B$. The cost to solve
the $j$-th sketched subproblem is dominated 
by the downsampled MTTKRP operation 
$\textrm{mat}(\scr T, j) S^\top S U_{\neq j}$, which has
runtime
\begin{equation}
O\paren{\nnz{\textrm{mat}(\scr T, j) S^\top} N R}.
\label{eq:ds_mttkrp_runtime}
\end{equation}
As Figure \ref{fig:mttkrp_definition} (bottom) illustrates, 
$\textrm{mat}(\scr T, j) S^\top$ typically has far 
fewer nonzeros than $\scr T$, enabling sampling to 
reduce the computation cost in 
Equation \eqref{eq:mttkrp_runtime}. To select indices
to sample, we implement two algorithms 
that involve the \emph{leverage scores} of the design matrix
\cite{cheng_spals_2016, larsen_practical_2022, bharadwaj2023fast}. 
Given a matrix $A \in \RR^{I \times R}$, the leverage score of row $i$ is given by 
\begin{equation}
\ell_i = A\br{i,:} (A^\top A)^+ A \br{i,:}^\top.
\label{eq:lev_definition}
\end{equation}

These scores induce a probability distribution over the rows of matrix $A$, which we can interpret as a measure of importance. As the following theorem from Larsen and Kolda \cite{larsen_practical_2022} 
(building on similar results by Mahoney and Drineas
\cite{drineas_mahoney_rnla_survey}) shows, 
sampling from either 
the exact or approximate distribution of statistical leverage
guarantees, with high probability, that the solution
to the downsampled problem has low residual with respect
to the original problem.
\begin{theorem}[Larsen 
and Kolda \cite{larsen_practical_2022}]
Let $S \in \RR^{J \times I}$ be a sampling matrix 
for $A \in \RR^{I \times R}$ 
where each row $i$ is sampled
i.i.d. with probability $p_i$. Let 
$\beta = \min_{i \in \br{I}} (p_i R / \ell_i)$. For a constant $C$ and any 
$\varepsilon, \delta \in (0, 1)$, let the sample count be
\[
J = \frac{R}{\beta} \max \paren{ C \log \frac{R}{\delta}, 
\frac{1}{\varepsilon \delta}}.
\]
Letting $\tilde X = \textrm{argmin}_{\tilde X} \norm{SA \tilde X - SB}_F$, 
we have 
\[
\norm{A \tilde X - B}^2_F 
\leq (1 + \varepsilon) \min_X \norm{A X - B}_F^2.
\]
with probability at least $1 - \delta$.
\label{thm:leverage_scores}
\end{theorem}
Here, $\beta \leq 1$ quantifies deviation of the
sampling probabilities from the exact leverage 
score distribution, with a higher sample
count $J$ required as the deviation increases. The STS-CP samples
from the exact leverage distribution with $\beta = 1$, achieving
higher accuracy at the expense of increased sampling time.
CP-ARLS-LEV samples from an approximate distribution with 
$\beta < 1$.

Sketching methods for tensor decomposition
have been extensively investigated \cite{cheng_spals_2016,ahle_treesketch,
larsen_practical_2022, malik_more_2022, 
practical_randomized_cp}, both in theory 
and practice. Provided an appropriate sketch 
row count $J$ and assumptions common in the optimization
literature, rigorous convergence guarantees for randomized 
ALS can be derived \cite{gittens_adaptive_sketching}.

\section{Related Work}

\subsection{High-Performance ALS CP Decomposition}
Significant effort has been devoted to 
optimizing the shared-memory MTTKRP 
using new data structures for 
the sparse tensor, cache-blocked
computation, loop reordering strategies, and methods that 
minimize data races between threads \cite{smith_csf, 
nisa_csf_mm, parti, oom_sparse_mttkrp,
phipps2019software,
wijeratne2023dynasor, kanakagiri2023minimum}. 
Likewise, several works provide high-performance algorithms for ALS CP decomposition in a distributed-memory setting.
Smith and Karypis provide an algorithm that distributes load-balanced chunks of the sparse tensor to processors
in an $N$-dimensional Cartesian topology \cite{smith_medium-grained_2016}. 
Factor matrices are shared among slices of the topology that require
them, and each processor computes a local MTTKRP before reducing results with a subset of 
processors. The SPLATT library \cite{smith_splatt_2015} implements this communication strategy and 
uses the compressed sparse fiber (CSF) 
format to accelerate local sparse MTTKRP computations 
on each processor.

Ballard et al. \cite{ballard_parallel_2018}. 
use a similar communication strategy to compute the MTTKRP involved 
in dense nonnegative CP decomposition. 
They further introduce a dimension-tree algorithm that reuses
partially computed terms of the MTTKRP between ALS optimization problems. DFacTo \cite{choi_dfacto_2014}
instead reformulates the MTTKRP as a sequence of sparse matrix-vector products (SpMV), 
taking advantage of extensive research optimizing the SpMV kernel. Smith and Karypis \cite{smith_medium-grained_2016}
note, however,
that DFacTo exhibits significant communication overhead. Furthermore,
the sequence of SpMV operations cannot take advantage of access locality within rows of the dense factor
matrices, leading to more cache misses than strategies based on sparse-matrix-times-dense-matrix-multiplication (SpMM). GigaTensor \cite{kang_gigatensor_2012} uses the MapReduce 
model in Hadoop to scale to distributed, fault-tolerant clusters. Ma and Solomonik \cite{ma_efficient_2021} use pairwise
perturbation to accelerate CP-ALS, reducing the cost of MTTKRP computations when ALS is sufficiently close
to convergence using information from prior rounds.

Our work investigates variants of the Cartesian data distribution scheme adapted for a
downsampled MTTKRP. We face challenges adapting either specialized data 
structures for the sparse tensor or dimension-tree algorithms. By extracting arbitrary nonzero elements 
from the sparse tensor, randomized sampling destroys the advantage conferred by 
formats such as CSF. Further, each least-squares solve requires a fresh set of rows
drawn from the Khatri-Rao product design matrix, which prevents efficient reuse of results 
from prior MTTKRP computations.

Libraries such as the Cyclops Tensor Framework (CTF) \cite{solomonik_massively_2014} automatically parallelize
distributed-memory contractions of both sparse and dense tensors. SpDISTAL \cite{yadav_spdistal_2022} 
proposes a flexible domain-specific language to schedule 
sparse tensor linear algebra on a cluster, including the MTTKRP operation. The randomized 
algorithms investigated here could be implemented on 
top of either library, but it is unlikely that current tensor 
algebra compilers can automatically produce 
the distributed samplers and optimized communication schedules that
we contribute.

\subsection{Alternate Sketching Algorithms and Tensor Decomposition
Methods}
\label{sec:alternate_sketches}
Besides leverage score sampling, popular
options for sketching Khatri-Rao products
include Fast Fourier Transform-based 
sampling matrices \cite{jin_fjlt} and structured random sparse matrices (e.g. Countsketch) \cite{ahle_treesketch, diao_kronecker_sketch}. 
The former method, however, introduces fill-in when applied to the sparse matricized tensor $\textrm{mat}(\scr T, j)$. Because the runtime of the downsampled MTTKRP 
is linearly proportional to the nonzero count of
$\textrm{mat}(\scr T, j) S^\top$, the advantages of 
sketching are 
lost due to fill-in. 
While Countsketch operators do
not introduce fill, they still require access to all 
nonzeros of the sparse
tensor at every iteration, which is expensive
when $\nnz{\scr T}$ ranges from hundreds of millions to
billions.

Other algorithms besides ALS exist for large sparse tensor decomposition. Stochastic gradient descent
(SGD, investigated by Kolda and Hong \cite{kolda_stochastic_2020}) 
iteratively improves CP factor matrices 
by sampling minibatches of indices from $\scr T$, computing 
the gradient of a loss function at those indices with respect to the factor matrices, and adding a step in the direction of the gradient to the factors. Gradient methods are flexible enough 
to minimize a variety of loss functions besides the Frobenius
norm error \cite{gcp_kolda}, but 
require tuning additional parameters (batch size, learning rate) and a distinct parallelization strategy. 

%The CCD++ algorithm \cite{yu_scalable_2012} 
%extended to tensors keeps all but one rank-1 component of the decomposition fixed and optimizes for the remaining 
%component, in contrast to ALS which keeps all but one factor matrix fixed. 

\section{Distributed-Randomized CP Decomposition}
\label{sec:dist_cp_decomposition}
\begin{figure*}
    \centering
    \includegraphics[trim={0.8cm 20.7cm 0cm 1.5cm},scale=0.9,clip]{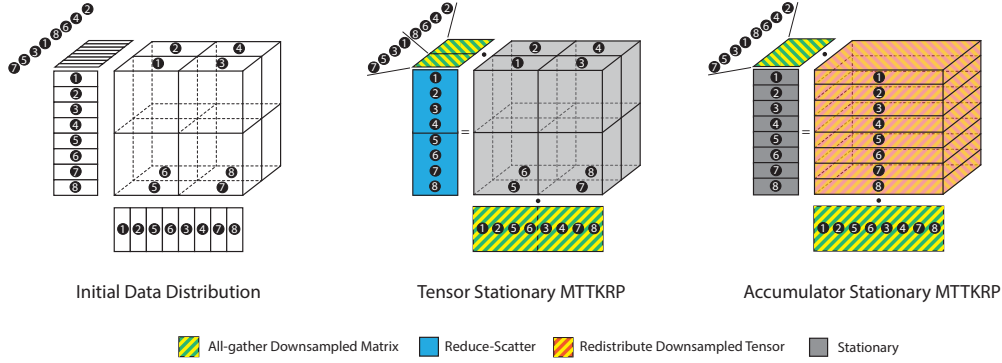}
    \vspace{-20pt}
    \caption{Initial data distribution and downsampled MTTKRP data movement for a 3D tensor, $P=8$ processors. Rectangles along each side of the tensor illustrate 
    factor matrices corresponding to each mode, divided by block rows among
    processors. Each black circle denotes the processor owning a block of a matrix or 
    tensor; multiple circles on an object 
    indicate replication of a piece of data. Colors / shading indicate communication collectives.} 
   \label{fig:mttkrp_strategies}
   \Description[Initial Data Distribution and
   Communication Strategies]{Three cubes, each partitioned
   according in a 3D grid with three rectangles aligned
   with each side. Cells of the cube and rectangles are
   labeled with processors that own each piece of data.
   Left cube: the initial data distribution. Middle cube:
   the tensor is not communicated, while all three of
   the factors are communicated. Right: Only two of
   the factor matrices and potentially the tensor are
   communicated, while one factor stays in place.}
   \vspace{-7pt}
\end{figure*}
In this section, we distribute Algorithm 
\ref{alg:als_algorithm} to $P$ processors when 
random sampling is used to solve the 
least-squares problem on line 6. Figure 
\ref{fig:mttkrp_strategies} (left) shows the initial
data distribution of our factor matrices and tensor
to processors, which are arranged in a hypercube
of dimensions $P_1 \times ... \times P_N$ with
$\prod_i P_i = P$. Matrices $U_1, ..., U_N$ 
are distributed by block rows 
among the processors to ensure an even division of
computation, and we denote by $U_i^{(p_j)}$ the
block row of $U_i$ owned by processor $p_j \in \br{P}$. 
We impose that all processors can access the Gram matrix $G_i$ of each factor $U_i$, which is computed 
by an \verb|Allreduce| of the $R \times R$ 
matrices $U_i^{(p_j)\top} U_i^{(p_j)}$ across 
$p_j \in \br{1, ..., P}$. Using these matrices,
the processors
redundantly compute the overall Gram matrix $G$
through Equation \eqref{eq:efficient_G_formula}, and
by extension $G^+$.

With these preliminaries, each processor takes
the following actions to execute steps 6-8 of
Algorithm \ref{alg:als_algorithm}:

\begin{enumerate}
    \item \textbf{Sampling and All-gather: } 
    Sample rows of $U_{\neq j}$ 
    according to the leverage-score distribution and
    \verb|Allgather| the rows to processors who 
    require them. For non-randomized ALS, no sampling
    is required.
 
    \item \textbf{Local Computation: } Extract the corresponding nonzeros from the local
    tensor owned by each processor and execute the
    downsampled MTTKRP, a sparse-dense matrix multiplication.

    \item \textbf{Reduction and Postprocessing: } Reduce the accumulator of the sparse-dense
    matrix multiplication across processors, if
    necessary, and post-process
    the local factor matrix slice by multiplying with 
    $G^+$. Renormalize the factor matrix columns and update 
    sampling data structures.
\end{enumerate}
Multiple prior works establish the correctness of this schedule \cite{smith_medium-grained_2016,ballard_parallel_2018}. We now examine
strategies for drawing samples (step 1),
communicating factor matrix rows (steps 2 and 3), and
performing local computation efficiently (step 2)
tailored to the case of randomized least-squares.

\subsection{New Distributed Sampling Strategies}
\label{sec:new_sampling_strategies}
Table \ref{tab:sampler_compute_comm} 
gives the asymptotic per-processor computation 
and communication
costs to draw $J$ samples in our 
distributed versions of CP-ARLS-LEV and STS-CP.
We give detailed descriptions, as well as
pseudo-code, for each sampling strategy 
in Appendices
\ref{sec:cp_arls_lev_sampling} and
\ref{sec:sts_cp_sampling}. 
In this section, we briefly describe the accuracy characteristics 
and communication / computation patterns for 
each method. Table \ref{tab:sampler_compute_comm}
does not include the costs to construct the sampling data structures in each algorithm, 
which are subsumed asymptotically 
by the matrix-multiplication 
$U_i^{(p_j)} \cdot G^+$ on each processor (step 3). The costs 
of all communication collectives are taken from Chan et al.\
\cite{chan_collective_2007}. 
\newline

\begin{table}
\begin{tabular}{@{}llll@{}}
\toprule
\textbf{Sampler} & \textbf{Compute}  & \textbf{Messages} & \textbf{Words Sent/Recv} \\ \midrule
d-CP-ARLS-LEV      & $JN / P$          & $P$                    & $JN / P$                        \\
d-STS-CP           & $(JN / P) R^2 \log P$ & $NP \log P$             & $(J / P) NR \log P$             \\ \bottomrule
\end{tabular}
\vspace{5pt}
\caption{Asymptotic Per-Processor Costs to Draw $J$
Samples}
\label{tab:sampler_compute_comm}
\end{table}

\begin{table}
\centering
\begin{tabular}{@{}llll@{}}
\toprule
\textbf{Schedule} & \textbf{Words Communicated / Round} \\
\midrule
Non-Randomized TS & $2NR\paren{\textstyle\prod_{k=1}^N I_k/P}^{1/N}$\\ 
Sampled TS & $NR\paren{\textstyle\prod_{k=1}^N I_k/P}^{1/N}$
\\ 
Sampled AS & $JRN(N-1)$ \\ 
\bottomrule
\end{tabular}
\vspace{5pt}
\caption{Communication Costs for Downsampled MTTKRP}
\label{tab:mttkrp_comm}
\end{table}

\textbf{CP-ARLS-LEV: } The CP-ARLS-LEV algorithm
by Larsen and Kolda \cite{larsen_practical_2022}
\textit{approximates} the leverage scores in Equation
\eqref{eq:lev_definition} by the product of leverage 
scores for each factor matrix $U_1, ..., U_N$. The
leverage scores of the block row $U_i^{(p_j)}$ owned
by processor $p_j$ are approximated by
\[
\tilde \ell^{(p_j)} = \textrm{diag}\paren{U_i^{(p_j)}G_i^+ U_i^{(p_j)\top}}
\]
which, given the replication of $G_i^+$,
can be constructed independently by each processor in
time $O\paren{R^2 I_i / P}$. The resulting 
probability vector, which is distributed among $P$ 
processors, can be sampled in expected time $O(J / P)$,
assuming that the sum of leverage scores distributed
to each processor is roughly equal (see Section
\ref{sec:load_balance} on load balancing for methods 
to achieve this). Multiplying by $(N-1)$ to sample independently 
from each matrix held constant, we get an asymptotic computation cost $O(JN / P)$ for the 
sampling phase. Processors
exchange only a constant multiple of $P$ 
words to communicate the sum of leverage
scores that they hold locally and the exact 
number of samples
they must draw, as well as a cost $O(JN / P)$ to
evenly redistribute / postprocess the final sample
matrix. While this algorithm is computationally efficient,
it requires $J = \tilde O(R^{N-1} / (\varepsilon \delta))$ 
to achieve the 
$(\epsilon, \delta)$-guarantee from Theorem
\ref{thm:leverage_scores}, which may lead to a higher
runtime in the distributed-memory MTTKRP. Larsen and Kolda note 
that CP-ARLS-LEV sampling can be implemented
without any communication at all if an entire factor matrix 
is assigned uniquely to a single processor, which can compute leverage 
scores and draw samples independently \cite{larsen_practical_2022}.
That said, assigning an entire factor matrix to a single processor
incurs higher communication costs in the MTTKRP phase of the 
algorithm and may be infeasible under tight memory constraints,
leading to our adoption of a block-row distribution for the factors. 
\newline

\begin{figure}
    \centering
    \begin{tikzpicture}[scale=0.55, level/.style={sibling distance=50mm/#1}]
    \node [very thick, circle, draw] (r){}
        child {node [circle, draw] (v1){}
            child {node [circle, draw] (v3){}  
                child {node [circle, draw] (v4){}} 
                child {node [circle, draw] (v5){}}
            } 
            child { node [circle, draw] (v6){}
                child {node [circle, draw] (v7){}} 
                child {node [circle, draw] (v8){}}
            }
        } 
        child [->, very thick] {node [circle, draw] (v9){} 
            child {node [circle, draw] (v10){}  
                child [-, thin] {node [circle, draw] (v11){}} 
                child {node [circle, draw] (v12){}}
            }
            child [-, thin] {node [circle, draw] (v13){} 
                child {node [circle, draw] (v14){}} 
                child {node [circle, draw] (v15){}}
            } 
        };
        \node[anchor=east, xshift=-15pt] at (v3.west) {$L$};
        \draw [dashed, shorten >= -90pt, shorten <=-25pt] (v3) -- (v13);

        \node[anchor=south] at (r.north) (l1) {$p_1, p_2, p_3, p_4$}; 
        
        \node[anchor=south east] at (v1.north west) (l2) {$p_1, p_2$}; 
        \node[anchor=south west] at (v9.north east) {$p_3, p_4$}; 
        
        \node[anchor=south east] at (v3.north west)  (l3) {$p_1$}; 
        \node[anchor=south west] at (v6.north east)  {$p_2$};  
        \node[anchor=south east] at (v10.north west)  {$p_3$};  
        \node[anchor=south west] at (v13.north east) {$p_4$};

        \node[anchor=center, xshift=13em, yshift=3.5em, text width=1.5cm, align=center] at (r.west) (owner) {\footnotesize \textbf{Sample Owner}};
        \path let \p1 = (l1), \p2 =(owner) in node at (\x2, \y1) (rw1) {\footnotesize Step 1: $p_2$};  
        \path let \p1 = (l2), \p2 =(owner) in node at (\x2, \y1) (rw2) {
        \footnotesize Step 2: $p_3$}; 
        \path let \p1 = (l3), \p2 =(owner) in node at (\x2, \y1) (rw3) {
        \footnotesize Step 3: $p_3$};
        \path let \p1 = (l3), \p2 =(owner) in node at (\x2, \y1-4.5em) (rw4) {
        \footnotesize Local Compute};
        
        \draw[->] (rw1) -> (rw2); 
        \draw[->] (rw2) -> (rw3); 
        \draw[->] (rw3) -> (rw4);

  \tikzset{
    mynode/.style={
      rectangle,
      draw,
      minimum width=1.38cm,
      minimum height=0.5cm,
      on chain,
    }
  }

  \begin{scope}[start chain=going right, node distance=0pt]
  \node [anchor=north west, xshift=-4pt, yshift=-5pt] at (v4.south west) [mynode] {\footnotesize $U_1^{(p_1)}$};
    \foreach \x in {2,3,4} {
      \node [mynode] {\footnotesize $U_1^{(p_\x)}$};
    }
  \end{scope}
\end{tikzpicture}
    \caption{Example random walk in STS-CP to draw a single sample index from 
    matrix $U_1$, distributed to $P=4$ processors. Annotations on the 
    tree (left) indicate processors that share data for each node. 
    The schedule to the right
    indicates the processor that owns the sample at each stage of the
    random walk. The sample begins randomly at $p_2$, then branches left to
    $p_3$ ($p_4$ shares node data and could also have been selected), involving
    communication of a vector corresponding to the sample from $p_2$ to $p_3$. The sample 
    remains at $p_3$ for the remainder of the walk.}
    \label{fig:sts_cp_parallelization}
    \Description[A Distributed Random Walk Down a Binary 
    Tree]{A tree with a single path highlighted showing
    the path that a random walk takes between 
    processors. Each node is annotated with a processor
    that owns it. Below the tree, a partitioned rectangle
    shows the factor matrix blocks owned by each processor.}
\end{figure}
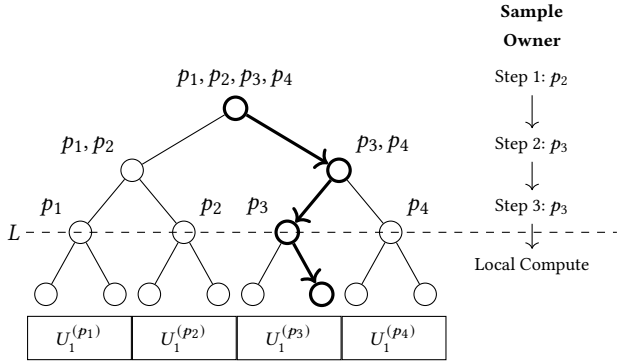

\textbf{STS-CP: } The STS-CP algorithm 
\cite{bharadwaj2023fast} samples from the exact 
leverage distribution by executing a random walk on
a binary tree data structure once for each of the $N-1$ factor
matrices held constant. Each leaf of the binary tree
corresponds to a block of $R$ rows from a factor matrix
$U_i$ and holds the $R \times R$ Gram matrix of that
block row. Each internal node $v$ holds a matrix 
$G^v$ that is the sum of the matrices held 
by its children. Each sample begins with a unique
vector $h$ at the root of the tree. At each non-leaf node $v$, the algorithm computes
$(h^\top G^{L(v)} h) / (h^\top G^{v} h)$. If this quantity is greater than
a random number $r$ unique to each sample, the algorithm sends the sample to the left
subtree, and otherwise the right subtree. The process repeats until the random walk 
reaches a leaf and a row index is selected.

We distribute the data structure and the random walk as shown in Figure 
\ref{fig:sts_cp_parallelization}. 
We assume that $P$ is a power of two to simplify our
description, but our implementation makes no such 
restriction. Each processor $p_j$ owns a 
subtree of the larger tree that corresponds to 
their block row $U_i^{(p_j)}$.
The roots of these subtrees all occur at the same depth 
$L = \log P$. Above level $L$,
each node stores $2 \log P$ additional matrices,
$G^v$ and $G^{L(v)}$, for each ancestor node $v$
of its subtree. 

To execute the random walks, each sample is 
assigned randomly to a processor which evaluates the branching threshold at the 
tree root. Based on the direction of the branch, the sample 
and corresponding vector $h$ are routed to 
a processor that owns the required node information, and the process repeats until the 
walk reaches level $L$. The remaining steps do not
require communication.

The replication of node information above level $L$ requires communication 
overhead $O(R^2 \log P)$ using the classic
bi-directional exchange algorithm for \verb|Allreduce| \cite{chan_collective_2007}. For
a batch of $J$ samples, each level of the tree requires $O(JR^2)$ FLOPs to evaluate
the branching conditions. Under the assumption that the final sampled rows 
are distributed evenly to processors, the computation and communication at
each level are load balanced in expectation. Each processor has 
expected computation cost $O((J / P) N R^2 \log P)$ over all levels of the tree and all matrices $U_i, 1 \leq i \leq N, i \neq k$. 
Communication of samples between tree levels is accomplished 
through \verb|All-to-allv| collective calls, requiring
$O(NP \log P)$ messages and $O((J / P) N R \log P)$ words sent / received in expectation
by each processor. 

\subsection{A Randomization-Tailored MTTKRP Schedule}
\label{sec:communication_schedules}
The goal of this section is to demonstrate that
an optimal communication schedule for non-randomized ALS may
incur unnecessary overhead for the randomized 
algorithm. In response, we will use a schedule 
where all communication costs scale with the number
of random samples taken, enabling 
the randomized algorithm to decrease
\textit{communication costs} as well as computation. Table \ref{tab:mttkrp_comm}
gives lower bounds on the communication
required for each schedule we consider, and
we derive the exact costs in
this section.

The two schedules that we consider
are ``tensor-stationary", where factor matrix rows 
are gathered and reduced across a grid, and 
``accumulator-stationary", where no reduction takes place. 
These distributions were compared by Smith and Karypis 
\cite{smith_medium-grained_2016} under the names 
``medium-grained" and ``course-grained", respectively. Both
distributions exhibit, under an even distribution  
of tensor nonzero entries and leverage scores to 
processors, ideal expected computation scaling. 
Therefore, we focus our analysis on communication.
We begin by deriving the communication costs for non-randomized
ALS under the tensor-stationary communication schedule, which
we will then adapt to the randomized case.

Although our input tensor is sparse, we model the \textbf{worst-case}
communication costs for the dense factor matrices 
with standard \verb|Allgather|
and \verb|Reduce-scatter| primitives. For non-randomized (exact) ALS, 
the cost we derive matches that given by Smith and Karypis in their 
sparse tensor decomposition
work \cite{smith_medium-grained_2016}. 
Furthermore consider the extremely sparse Reddit tensor, (nonzero fraction $4 \times 10^{-10}$ \cite{smith_frostt_2017}),
which nonetheless exhibits an average of 571 nonzeros per fiber along 
the longest tensor mode and an average of 26,000 nonzeros per fiber 
aligned with the shortest tensor mode. The high 
per-fiber nonzero count induces a practical
communication cost comparable to the 
worst-case bounds, a feature that Reddit shares 
with other datasets in Table \ref{tab:datasets}.
\newline

\textbf{Exact Tensor-Stationary:}
\label{section:tensor_stationary}
The tensor-stationary MTTKRP algorithm is communication-optimal for dense CP
decomposition \cite{ballard_parallel_2018} and outperforms several other methods 
in practice for non-randomized sparse CP decomposition.
\cite{smith_medium-grained_2016}. The middle image of Figure \ref{fig:mttkrp_strategies} illustrates the approach.
During the $k$-th optimization problem in a round of ALS, each processor does the following: 
\begin{enumerate}
    \item For any $i \neq k$, participates in an \verb|Allgather| of all blocks $U_i^{(p_j)}$ for all
    processors $p_j$ in a slice of the 
    processor grid aligned
    with mode $k$. 

    \item Executes an MTTKRP with locally owned
    nonzeros and the gathered row blocks. 

    \item Executes a \verb|Reduce-scatter| with the MTTKRP result along a slice of the
    processor grid aligned with mode $j$, storing
    the result in $U_k^{(p_j)}$
\end{enumerate}
For non-randomized ALS, the gather step
must only be executed once per round and can be cached. 
Then the communication
cost for the \verb|All-gather| and \verb|Reduce-scatter|
collectives summed over all $k=1...N$ is
\[
    2\textstyle\sum_{k=1}^N I_k R/P_k.
\]
To choose the optimal grid dimensions $P_k$, we minimize the expression above subject to the constraint $\prod_{k=1}^N P_k = P$. Straightforward application
of Lagrange multipliers leads to the optimal grid dimensions 
\[
    P_k = I_k\paren{ P /\textstyle\prod_{i=1}^N I_i}^{1/N}.
\]
These are the same optimal grid dimensions reported 
by Ballard et al.\ \cite{ballard_parallel_2018}. The communication under this optimal grid is 
\[
    2NR\paren{\textstyle\prod_{k=1}^N I_k/P}^{1/N}.
\]

\textbf{Downsampled Tensor-Stationary}:
As Figure \ref{fig:mttkrp_strategies} illustrates,
only factor matrix rows that are selected by the
random sampling algorithm need to be gathered by
each processor in randomized CP decomposition. 
Under the assumption that 
sampled rows are evenly distributed among 
the processors, the 
expected cost of gathering rows reduces to 
$JR(N - 1) / P_k$ within slices along mode $k$. 
The updated communication cost under the optimal grid
dimensions derived previously is 
$$\frac{R\paren{\textstyle\prod_{k=1}^N I_k}^{1/N}}{P^{1/N}} \br{N+ \sum_{k = 1}^N 
\frac{J(N-1)}{I_k}}.$$
The second term in the bracket arises from \verb|Allgather| collectives
of sampled rows, which is small if $J \ll I_k$ for all $1 \leq k \leq N$. 
The first term in the bracket arises from the \verb|Reduce-scatter|, which is unchanged by the
sampling procedure. Ignoring the second term in the
expression above gives the second entry of Table \ref{tab:mttkrp_comm}. 

Observe that this randomized method spends the same time 
on the reduction as the non-randomized schedule 
while performing significantly less computation, leading to diminished arithmetic intensity. On the other hand, this distribution may be
optimal when the tensor dimensions $I_k$ are small or
the sample count $J$ is high enough. 
\newline

\textbf{Downsampled Accumulator-Stationary}: As shown by Smith and
Karypis \cite{smith_medium-grained_2016}, the accumulator-stationary data
distribution performs poorly for non-randomized ALS. In the worst case, each 
processor requires access to all entries from all factors $U_1, ..., U_N$, 
leading to high communication and memory overheads. On the other hand, we demonstrate
that this schedule may be optimal for \textit{randomized} ALS on tensors where the
sample count $J$ is much smaller than the
tensor dimensions. The rightmost image in 
Figure \ref{fig:mttkrp_strategies}
illustrates the approach, which avoids the 
expensive \verb|Reduce-scatter| collective. To optimize $U_k$, we 
keep the destination buffer for a block row of $U_k$ stationary on each 
processor while communicating only
sampled factor matrix rows and nonzeros of $\scr T$. Under this distribution, all sampled factor matrix rows must be gathered to all
processors. The cost of the gather step for a single round becomes 
$O(JRN(N-1))$ (for each of $N$ least-squares problems, we 
gather at most $J(N-1)$ rows of length $R$). 
Letting $S_1, ..., S_N$ be the 
sampling matrices for each ALS subproblem  
in a round, the number of nonzeros selected
in problem $j$ is 
$\textrm{nnz}(\textrm{mat}(\scr T, j) S_j)$.
These selected (row, column, value) triples 
must be redistributed
as shown in Figure \ref{fig:mttkrp_strategies}
via an \verb|All-to-allv| collective call.
Assuming that the source and destination for
each nonzero are distributed uniformly among
the processors, the expected cost of redistribution in
least-squares problem $j$ is 
$(3 / P) \textrm{nnz}(\textrm{mat}(\scr T, j) S_j^\top)$. 
The final communication cost is 
\begin{equation}    
JRN(N-1) + \frac{3}{P} 
\sum_{j=1}^N \textrm{nnz}(\textrm{mat}(\scr T, j) S_j^\top).
\label{eq:accum_stationary_cost}
\end{equation}
The number of nonzeros sampled varies 
from tensor to tensor even when the sample count
$J$ is constant. That said, the redistribution 
exhibits perfect scaling (in expectation) with the processor
count $P$. In practice, we avoid 
redistributing the tensor entries multiple times by storing $N$ different 
representations of the tensor aligned with each slice of the processor grid, 
a technique that competing packages (e.g. DFacto 
\cite{choi_dfacto_2014}, early 
versions of SPLATT \cite{smith_csf})
also employ. This optimization eliminates the 
second term in Equation \eqref{eq:accum_stationary_cost}, 
giving the communication cost in the third row of Table \ref{tab:mttkrp_comm}. More importantly, observe that all 
communication scales linearly with the sample count 
$J$, enabling sketching to improve \textit{both} 
the communication and computation efficiency of our
algorithm. On the other hand, the term $JRN(N-1)$ 
does not scale with $P$, and we expect that 
gathering rows becomes a communication bottleneck for high processor counts.

\subsection{Tensor Storage and Local MTTKRP}
\label{sec:local_tensor_storage}
As mentioned in Section \ref{sec:communication_schedules}, we store
different representations of the sparse tensor 
$\scr T$ 
across the processor
grid to decrease communication costs. Each corresponds 
to a distinct matricization 
$\textrm{mat}(\scr T, j)$ for $1 \leq j \leq N$ used in the 
MTTKRP (see Figure \ref{fig:mttkrp_definition}). For 
non-randomized ALS, a variety of alternate storage formats have 
been proposed to reduce the memory overhead and accelerate the
local computation. Smith and Karypis 
support a compressed sparse fiber format for the tensor in SPLATT
\cite{smith_csf, smith_splatt_2015}, and 
Nisa et al. \cite{nisa_csf_mm} propose a mixed-mode compressed
sparse fiber format as an improvement. These optimizations 
cannot improve the runtime of our randomized algorithms 
because they are not conducive to sampling random nonzeros from
$\scr T$.

\begin{figure}
    \centering
    \begin{tikzpicture}[scale=0.7,
mnd/.style n args={2}{transform shape, draw, 
minimum height=#1 cm, minimum width=#2 cm,anchor=center,font=\footnotesize}
]

\node[right, mnd={2}{4.8}, anchor=north west] at (0, 0) (spmat) {};

\definecolor{gray_cust}{RGB}{200,200,200} 
% Gray shading showing rows taken from the MTTKRP

\foreach \x in {1.0cm, 0.5cm, 1.8cm, 4.1cm}
{
\node[mnd={1.98}{0.2}, anchor=center, fill=gray_cust, draw=none,
xshift=\x] at (spmat.west) {};
}

% Arrows in the sparse tensor indicating the order in which nonzero
% entries are stored.
\draw[->, draw=darkgray] (spmat.north west) ++(0.20, - 0.05)  -- ++(0, -1.8); 
\draw[->, draw=darkgray] (spmat.north west) ++(0.50, - 0.05)  -- ++(0, -1.8); 
\draw[->, draw=darkgray] (spmat.north west) ++(0.80, - 0.05)  -- ++(0, -1.2); 

% Put points in the sparse tensor
\foreach \val in {(1.3, -0.5), (-1.6, -0.3), (-0.6, -0.5), 
(0.7, 0.3), (1.7, 0.5), (-1.9, 0.5)}
{
\draw[fill=black] (spmat.center) ++\val circle(0.08em);
}
\node[font=\footnotesize] at (spmat.center) {$\textrm{mat}(\scr T, 2)$};

% Second image 
\draw[->] (spmat.south) ++(0, -0.25)  -- ++(0, -1); 

\node[right, mnd={2}{2.5}, anchor=north, fill=gray_cust, yshift=-1.5cm] at (spmat.south) (spmat_reduced) {};

% Put arrows indicating the local storage order
\draw[->, draw=darkgray] (spmat_reduced.north west) ++(0.10, -0.16)  -- ++(2.3, 0); 
\draw[->, draw=darkgray] (spmat_reduced.north west) ++(0.10, -0.505)  -- ++(1.6, 0); 

% Put points in the sparse tensor
\foreach \val in {(-0.6, 0.5), (0.5, 0.5), (-0.1, -0.5)}
{
\draw[fill=black] (spmat_reduced.center) ++\val circle(0.08em);
}
\node[font=\footnotesize, yshift=-10] at (spmat_reduced.south) {$\textrm{mat}(\scr T, 2) S^\top$};
\node[right, anchor=north west, font=\Huge] at (spmat_reduced.north east) (cdot1) {$\cdot$};
\node[mnd={2.5}{1}, anchor=north west, fill=gray_cust] at (cdot1.north east) (U3) {};
\node[yshift=-0.5cm] at (U3.south) {\footnotesize $S (U_3 \odot U_1)$};

\node[left, anchor=north east, font=\large] at (spmat_reduced.north west) (coloneq) {$:=$};
\node[mnd={2}{1}, anchor=north east] at (coloneq.north west) (U2) {};
%\node[font=\footnotesize, yshift=-10] at (U2.south) 
%{$\tilde U_2$};

% Show the parallelization scheme 
\draw[draw=darkgray, dashed] (U2.north west) ++(0, -0.66)  -- ++(1.0, 0); 
\draw[draw=darkgray, dashed] (U2.north west) ++(0, -1.33)  -- ++(1.0, 0); 

\draw[draw=darkgray, dashed] (spmat_reduced.north west) ++(0, -0.66)  -- ++(2.5, 0); 
\draw[draw=darkgray, dashed] (spmat_reduced.north west) ++(0, -1.33)  -- ++(2.5, 0); 

% Label the thread blocks 
\node[font=\footnotesize, yshift=-6, xshift=-10] at (U2.north west) 
{$t_1$};
\node[font=\footnotesize, yshift=-20, xshift=-10] at (U2.north west) 
{$t_2$};
\node[font=\footnotesize, yshift=-34, xshift=-10] at (U2.north west) 
{$t_3$};

\end{tikzpicture}
    \caption{Shared-memory parallelization of 
    downsampled MTTKRP procedure. Nonzero 
    sparse coordinates in the sampled gray columns, 
    initially sorted by column, are selected and remapped 
    into a CSR matrix. The subsequent
    matrix multiplication is 
    parallelized to threads $t_1, t_2, t_3$ without atomic 
    operations or data races, since each thread is responsible 
    for a unique block of the
    output.}
    \label{fig:sparse_transpose}
    \Description[Sample Selection and Sparse Transpose]{
    A rectangle representing the matricized tensor with
    shaded columns selected by sampling. Below, three
    rectangles illustrate sparse-dense matrix multiplication
    with each thread assigned a block row of the output.
    Arrows inside the original matricized tensor and
    the downsampled version indicate the data ordering,
    which changes from column-major to row-major order
    as a result of transposition.}
\end{figure}
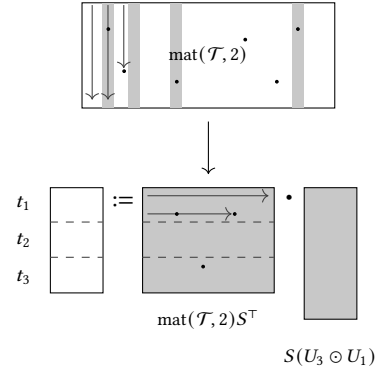

Instead, we adopt the approach shown in Figure 
\ref{fig:sparse_transpose}. The coordinates in each tensor 
matricization are stored in sorted order of their column indices,
an analogue of compressed-sparse-column (CSC) format. With this
representation, the random sampling 
algorithm efficiently selects columns of 
$\textrm{mat}(\scr T, j)$ corresponding to rows of the design matrix. 
The nonzeros in these columns are extracted and remapped to
a compressed sparse row (CSR) format through a ``sparse transpose"
operation. The resulting CSR matrix participates in 
the sparse-dense matrix multiplication, which
can be efficiently parallelized without data races 
on a team of shared-memory threads.

%The key to efficiency in the sparse matrix transpose is that
%the sampling process extracts only a small fraction of nonzero 
%entries from the entire tensor. We leave reducing 
%the memory footprint of 
%our randomized algorithms as future work.

\subsection{Load Balance}
\label{sec:load_balance}
To ensure load balance among processors, we randomly
permute the sparse tensor indices along each mode, a technique also used by SPLATT
\cite{smith_medium-grained_2016}. These permutations ensure that each 
processor holds, in expectation, an equal fraction of nonzero entries from
the tensor and an equal fraction of sampled nonzero entries. For
highly-structured sparse tensors, random 
permutations do not optimize processor-to-processor communication costs,
which packages such as Hypertensor \cite{hypertensor} minimize through 
hypergraph partitioning. As Smith and Karypis 
\cite{smith_medium-grained_2016} demonstrate empirically, 
hypergraph partitioning is slow and memory-intensive on large tensors.
Because our randomized implementations require just minutes on massive tensors
to produce decompositions comparable to non-randomized ALS,  
the overhead of partitioning outweighs the modest communication reduction it may produce. 

\section{Experiments}
\begin{table}
\centering
\begin{tabular}{lcll}
\toprule
Tensor & Dimensions & NNZ & Prep. \\ \midrule
Uber & $183 \times 24 \times 1.1K \times 1.7K$ & 3.3M & -\\
Amazon & $4.8M \times 1.8M \times 1.8M$ & 1.7B & -\\
Patents & $46 \times 239K \times 239K$ & 3.6B & - \\
Reddit & $8.2M \times 177K \times 8.1M$ & 4.7B & log \\
\bottomrule
\end{tabular}
\vspace{5pt}
\caption{Sparse Tensor Datasets}
\label{tab:datasets}
\end{table}
Experiments were conducted on CPU nodes of NERSC Perlmutter, a Cray HPE EX
supercomputer. Each node has 128 physical
cores divided between two AMD 
EPYC 7763 (Milan) CPUs. Nodes are linked by an 
HPE Slingshot 11 interconnect.

Our implementation is written in C++ and links with OpenBLAS 0.3.21 for
dense linear algebra. We use a simple Python wrapper around the C++
implementation to facilitate benchmarking. We use a hybrid 
of MPI message-passing and OpenMP shared-memory parallelism in
our implementation, which is available online
at \url{https://github.com/vbharadwaj-bk/rdist_tensor}.

Our primary baseline is the SPLATT , the Surprisingly Parallel 
Sparse Tensor Toolkit 
\cite{smith_medium-grained_2016, smith_splatt_2015}. 
SPLATT is a scalable CP decomposition package 
optimized for both communication costs
and local MTTKRP performance through innovative sparse
tensor storage structures. As a result, it remains one
of the strongest libraries for sparse tensor decomposition in 
head-to-head benchmarks against other libraries
\cite{rolinger_performance, nisa_csf_mm, kanakagiri2023minimum}. We 
used the default medium-grained
algorithm in SPLATT and adjusted the OpenMP
thread count for each tensor to achieve the
best possible performance to compare against.

Table \ref{tab:datasets} lists the sparse tensors used in our experiments, all sourced from 
the Formidable Repository of Open
Sparse Tensors and Tools (FROSTT) \cite{smith_frostt_2017}. 
Besides Uber, which was only used to verify accuracy
due to its small size, the Amazon, Patents, 
and Reddit tensors
are the only members of FROSTT at publication time 
with over 1 billion nonzero
entries. These tensors were identified to benefit the 
most from randomized sampling since the next 
largest tensor in the collection, NELL-1, has 12 times
fewer nonzeros than Amazon. We computed the logarithm
of all values in the Reddit tensor, consistent with
established practice \cite{larsen_practical_2022}.

\subsection{Correctness at Scale}
Table \ref{tab:accuracies} gives the average fits (5 trials) of decompositions produced by our 
distributed-memory algorithms. The fit \cite{larsen_practical_2022} between the decomposition 
$\tilde{\scr T} = \br{\sigma; U_1, ..., U_N}$ and the ground-truth
$\scr T$ is defined as 
\[
\textrm{fit}(\tilde{\scr T}, \scr T) = 
1 - \frac{\norm{\tilde{\scr T} - \scr T}_F}{\norm{\scr T}_F}.
\]
A fit of 1 indicates perfect agreement between the decomposition
and the input tensor. We used $J=2^{16}$ for our 
randomized algorithms to test our implementations on 
configurations identical to those in prior work 
\cite{larsen_practical_2022, bharadwaj2023fast}. To test both the distributed-memory message passing and 
shared-memory threading parts of our implementation, we used 
32 MPI ranks and 16 threads per rank 
across 4 CPU nodes. We report accuracy for the 
accumulator-stationary versions
of our algorithms and checked that the tensor-stationary
variants produced the same mean fits. The ``Exact"
column gives the fits generated by SPLATT. ALS was run for 40 rounds on all tensors except Reddit, for which
we used 80 rounds.

The accuracy of both d-CP-ARLS-LEV and d-STS-CP match
the shared-memory prototypes in the original
works \cite{larsen_practical_2022, bharadwaj2023fast}.
As theory predicts, the accuracy gap between 
d-CP-ARLS-LEV and d-STS-CP widens at higher rank. The 
fits of our methods improves by increasing the sample count $J$ at the expense of higher sampling and MTTKRP runtime.
\begin{table}
\centering
\begin{tabular}{llrr|r}
\toprule
Tensor & $R$ & d-CP-ARLS-LEV & d-STS-CP & Exact\\ \midrule
\multirow{3}{*}{Uber}   
& 25   & 0.187 & 0.189 & 0.190 \\ 
& 50   & 0.211 & 0.216 & 0.218 \\ 
& 75   & 0.218 & 0.230 & 0.232 \\ 
\midrule
\multirow{3}{*}{Amazon}   
& 25   & 0.338 & 0.340 & 0.340 \\ 
& 50   & 0.359 & 0.366 & 0.366 \\ 
& 75   & 0.368 & 0.381 & 0.382 \\ 
\midrule
\multirow{3}{*}{Patents}   
& 25   & 0.451 & 0.451 & 0.451 \\ 
& 50   & 0.467 & 0.467 & 0.467 \\ 
& 75   & 0.475 & 0.475 & 0.476 \\ 
\midrule
\multirow{3}{*}{Reddit}   
& 25   & 0.0583 & 0.0592 & 0.0596 \\ 
& 50   & 0.0746 & 0.0775 & 0.0783 \\ 
& 75   & 0.0848 & 0.0910 & 0.0922 \\ 
\bottomrule
\end{tabular}
\vspace{5pt}
\caption{Average Fits, $J=2^{16}$, 32 MPI Ranks, 4 Nodes}
\label{tab:accuracies}
\end{table}

\subsection{Speedup over Baselines}
\begin{figure}
    \centering
    \includegraphics[scale=0.40]{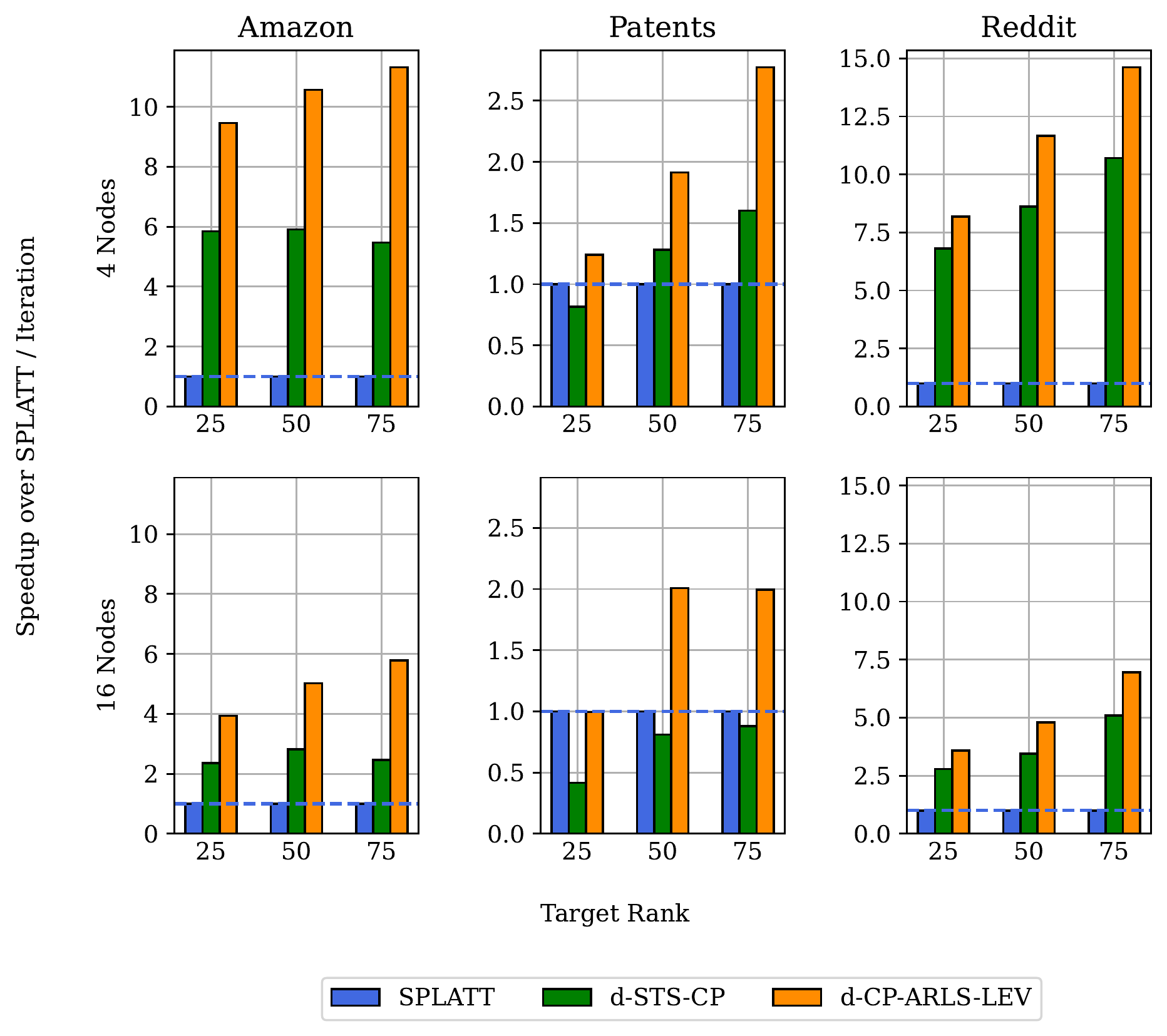}
    \caption{Average speedup per ALS iteration 
    of our distributed randomized 
    algorithms over SPLATT (5 trials, $J=2^{16}$).}
    \label{fig:splatt-speedup}
    \Description[Speedup of CP-ARLS-LEV and STS-CP over
    SPLATT]{Two rows of of colored bars showing the
    speedup of CP-ARLS-LEV and STS-CP 
    over SPLATT. Plots are shown for 4 nodes (top row)
    and 16 nodes (bottom row) on Amazon, Patents, and
    Reddit for ranks 25-75. The bar for SPLATT is 
    normalized to 1; speedups for STS-CP are at least 
    5x for Amazon and Reddit (maximum of 10x). CP-ARLS-LEV
    exhibits at least 9x speedup on Amazon and Reddit 
    (maximum of 15x). Speedups are smaller 
    (in the 1-2.5x range) for the denser Patents tensor.}
\end{figure}

Figure \ref{fig:splatt-speedup} shows the speedup of our randomized
distributed algorithm per ALS round over SPLATT at 4 nodes and 16 nodes. We used
the same configuration and sample
count for each tensor as Table
\ref{tab:accuracies}. On Amazon and Reddit at rank 25 and 4 nodes, 
d-STS-CP achieves a speedup in the range 
5.7x-6.8x while d-CP-ARLS-LEV achieves between 8.0-9.5x. We 
achieve our most dramatic speedup at rank 75 on the Reddit
tensor, with d-STS-CP achieving 10.7x speedup and d-CP-ARLS-LEV
achieving 14.6x. Our algorithms achieve less speedup compared to SPLATT on
the denser Patents tensor. Here, a 
larger number 
nonzeros are selected by randomized sampling,
with a significant computation bottleneck in
the step that extracts and reindexes the 
nonzeros from the tensor. The bottom half 
of Figure \ref{fig:splatt-speedup} shows that 
d-STS-CP maintains at least a 2x
speedup over SPLATT even at 16 nodes / 
2048 CPU cores on Amazon and Reddit, but
exhibits worse speedup on the Patents
tensor. Table \ref{tab:accuracies} 
quantifies the
accuracy sacrificed for the speedup, which can be changed by
adjusting the sample count at each least-squares solve. As Figure 
\ref{fig:reddit_fit_vs_time} shows, both of our randomized algorithms 
make faster progress than SPLATT, 
with d-STS-CP producing a
comparable rank-100 decomposition 
of the Reddit tensor in under two minutes.
\subsection{Comparison of Communication Schedules}

\begin{figure}
    \centering
    \includegraphics[scale=0.40]{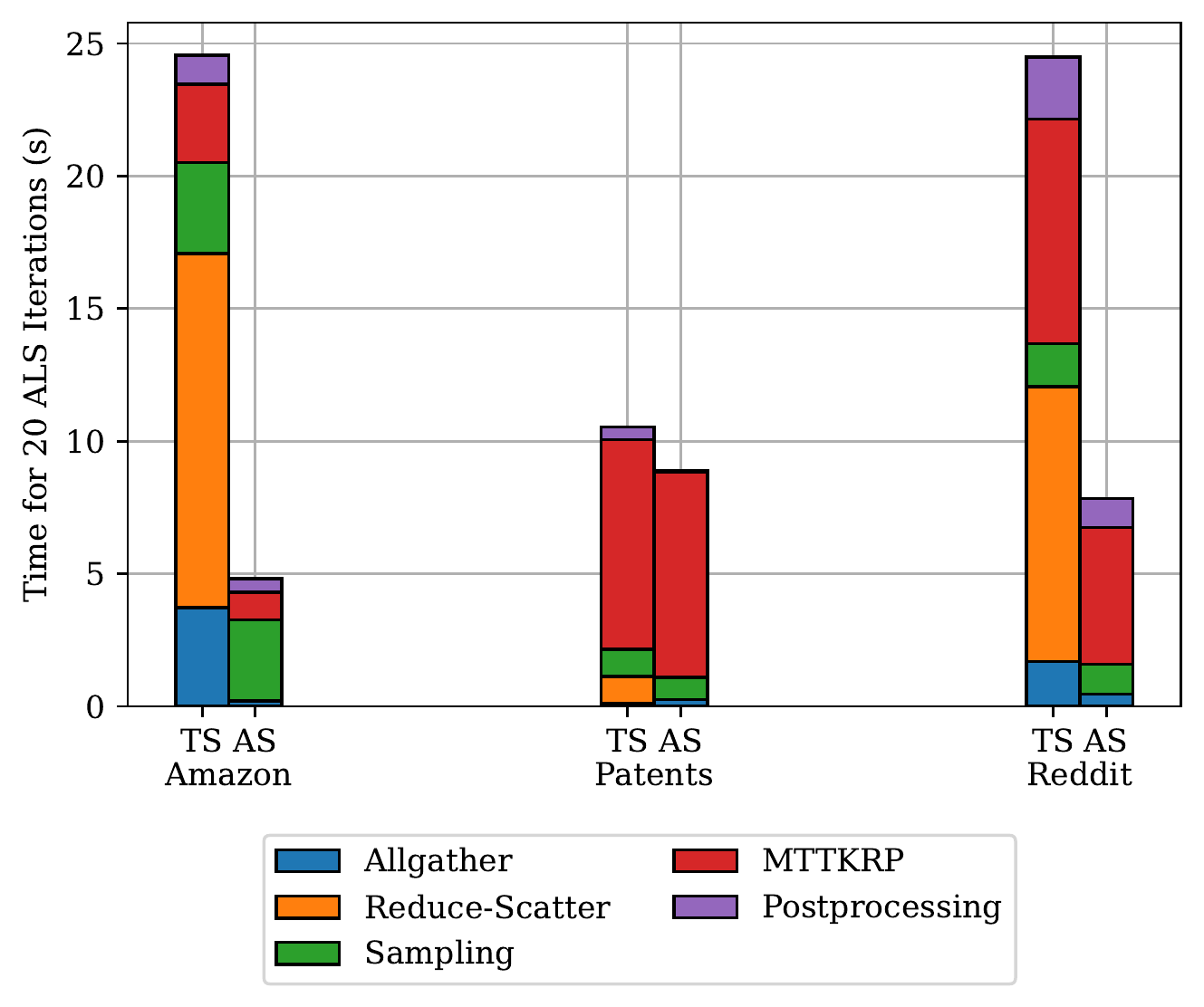}
    \caption{Average runtime (5 trials, $R=25$) 
    per activity
    for tensor-stationary and
    accumulator-stationary distributions
    with 32 MPI ranks over 4 nodes.}
    \label{fig:communication_comparison}
    \Description[Bar Charts Comparing Communication Schedules]{Three pairs of bars giving 
    the runtime for the tensor-stationary schedule
    and accumulator-stationary schedules on the
    three tensors. Each bar is broken down by activity.
    The runtime of reduction dominates the runtime for
    Amazon and Reddit, causing the tensor-stationary runtime
    to be much higher than accumulator-stationary. The
    bars are comparable in height for the Patents tensor.}
\end{figure}
Figure \ref{fig:communication_comparison}
breaks down the runtime per phase of the
d-STS-CP algorithm for the tensor-stationary 
and accumulator-stationary
schedules on 4 nodes. To illustrate the effect of sampling 
on the row gathering step, we gather all
rows (not just those sampled) for the tensor-stationary
distribution, a communication 
pattern identical to SPLATT. Observe that the 
\verb|Allgather| collective under the
accumulator-stationary schedule is significantly cheaper 
for Amazon and Reddit, since only sampled rows are
communicated. As predicted, the 
\verb|Reduce-scatter| collective
accounts for a significant fraction of
the runtime for the tensor-stationary
distribution on Amazon and Reddit, which have 
tensor dimensions in 
the millions. On both tensors, the runtime 
of this collective is greater than the time required
by all other phases combined in the 
accumulator-stationary schedule.
By contrast, both schedules perform comparably 
on Patents. Here, the
\verb|Reduce-scatter| cost is marginal
due to the smaller dimensions of the 
tensor. 

We conclude that sparse tensors with
large dimensions can benefit from the
accumulator-stationary distribution to
reduce communication costs, while the 
tensor-stationary distribution is
optimal for tensors with higher density
and smaller dimensions. The difference in MTTKRP 
runtime between the two schedules is further explored
in Section \ref{sec:load_balance_exps}.

\subsection{Strong Scaling and Runtime Breakdown} 
\label{sec:strong_scaling}
\begin{figure}
    \centering
    \includegraphics[scale=0.30]{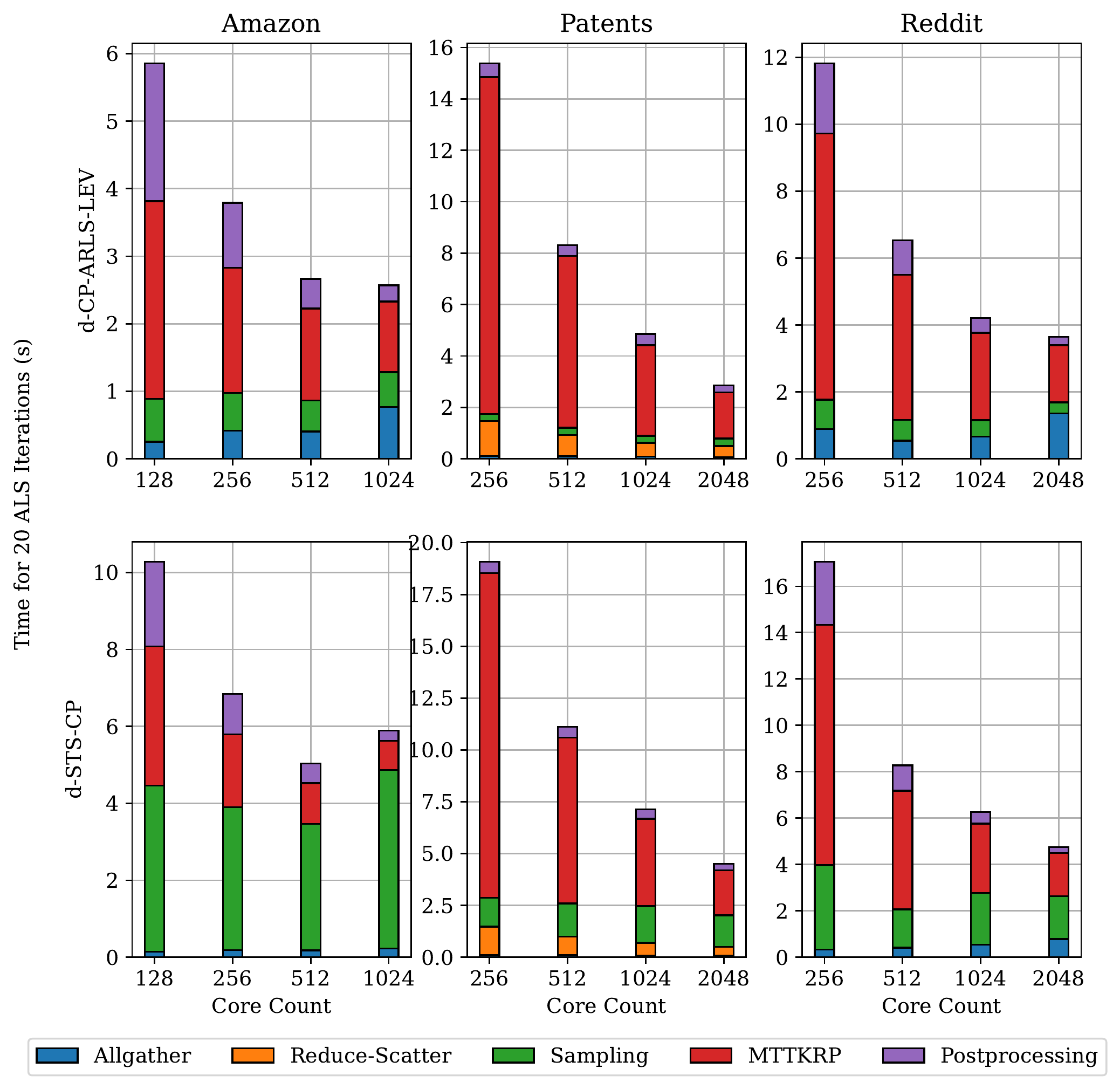}
    \caption{Average runtime (5 trials) per activity vs. CPU
    core count, $R=25$. Each node has 128 CPU cores, 
    and 8 MPI ranks were used per node.} 
    \label{fig:strong-scaling}
    \Description[Strong Scaling Bar Charts]{Two rows of
    three bar charts, each showing the runtime for increasing
    processor counts on the Amazon, Patents, and Reddit
    tensors. Scaling is near ideal for Patents and Reddit,
    while it slows by 1024 cores for Amazon. Bars are broken
    down by activity, which is dominated by the downsampled
    MTTKRP for Patents / Reddit. Sampling dominates the 
    runtime for STS-CP on Amazon, while all activities 
    require comparable time for CP-ARLS-LEV on Amazon.}
\end{figure}
Figure \ref{fig:strong-scaling} gives the runtime breakdown for our 
algorithms at varying core counts. Besides the \verb|All-gather| and
\verb|Reduce-scatter| collectives used to communicate rows of the factor
matrices, we benchmark time spent in each of the three phases identified
in Section \ref{sec:dist_cp_decomposition}: sample identification, execution of the downsampled MTTKRP, and post-processing factor matrices.

With its higher density, the Patents tensor has a significantly larger fraction 
of nonzeros randomly sampled at each linear least-squares solve. 
As a result, most ALS runtime is spent on the downsampled MTTKRP. The Reddit
and Amazon tensors, by contrast, spend a larger runtime portion 
on sampling and post-processing the factor matrices due to their larger
mode sizes. Scaling beyond 1024 cores 
for the Amazon tensor is impeded by the relatively high
sampling cost in d-STS-CP, a consequence of repeated
\verb|All-to-allv| collective calls. The high sampling 
cost is because the Amazon tensor has side-lengths 
in the millions along all tensor modes, leading to
deeper trees for the random walks in STS-CP.

\subsection{Weak Scaling with Target Rank} 
We measure weak scaling for our randomized
algorithms by recording the 
throughput (nonzero entries processed
in the MTTKRP per second of total algorithm runtime) as both the
processor count and target rank $R$ 
increase proportionally. We keep the ratio of
node count to rank $R$ constant at 16. 
We use a fixed sample count
$J=2^{16}$, and we benchmark the d-STS-CP 
algorithm to ensure minimal accuracy
loss as the rank increases. 

Although the FLOP count of the MTTKRP is linearly proportional to 
$R$ (see Equation \eqref{eq:mttkrp_runtime}), we expect the
efficiency of the MTTKRP to \textit{improve} with increased rank due
to spatial cache access locality in the longer factor matrix
rows, a well-documented phenomenon \cite{aktulga_tall_skinny}. 
On the other hand, the sampling runtime of the
d-STS-CP algorithm grows quadratically with the rank $R$ (see 
Table \ref{tab:sampler_compute_comm}). The net impact of these 
competing effects is determined by the density and
dimensions of the sparse tensor.
\begin{figure}
    \centering
    \includegraphics[scale=0.36]{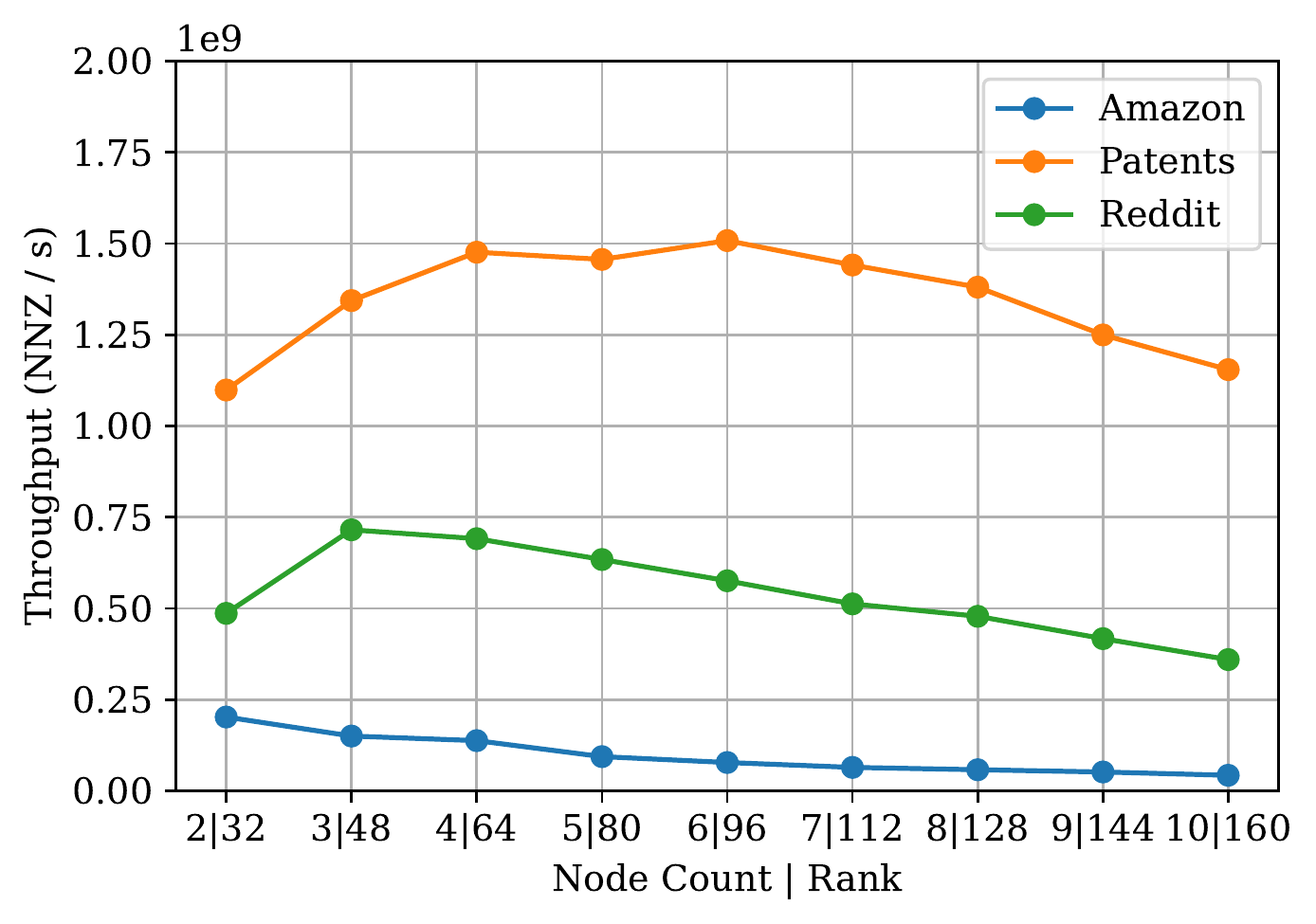}
    \caption{Average throughput (3 trials per
    data point) of the d-STS-CP algorithm 
    vs. increasing 
    node count and rank, measured
    as the average number of nonzeros iterated over in the MTTKRP
    per second of total algorithm runtime (higher is better). 
    Ideal scaling is a horizontal line. The ratio of node count to 
    rank was kept constant at 16. d-STS-CP 
    was chosen to preserve decomposition accuracy at high ranks.}
    \label{fig:weak_scaling}
    \Description[Weak Scaling Curves, Throughput vs. Node
    Count]{Three curves showing the throughput as both
    node count and rank increase from 2 nodes / rank-32 to
    10 nodes and rank-160. The throughput on Patents was
    highest (slightly above 1 billion nonzeros per second)
    followed by Reddit and Amazon (the latter just below
    250 million nonzeros per second). All three curves
    trend downwards at high node counts.}
\end{figure}

Figure \ref{fig:weak_scaling} shows the results of our weak scaling
experiments. Because ALS on 
the Amazon tensor spends 
a large fraction of time drawing samples
(see Figure \ref{fig:strong-scaling}), 
its throughput suffers with increasing
rank due to the quadratic 
cost of sampling. At the other extreme,
our algorithm spends little time sampling
from the Patents tensor with its
smaller dimensions, enabling throughput
to increase due to higher cache spatial
locality in the factor matrices. The experiments 
on Reddit follow a middle path between these extremes, 
with performance dropping slightly at high rank due to
the cost of sampling.

\subsection{Load Imbalance}
\label{sec:load_balance_exps}
\begin{figure}
    \centering
    \includegraphics[scale=0.48]{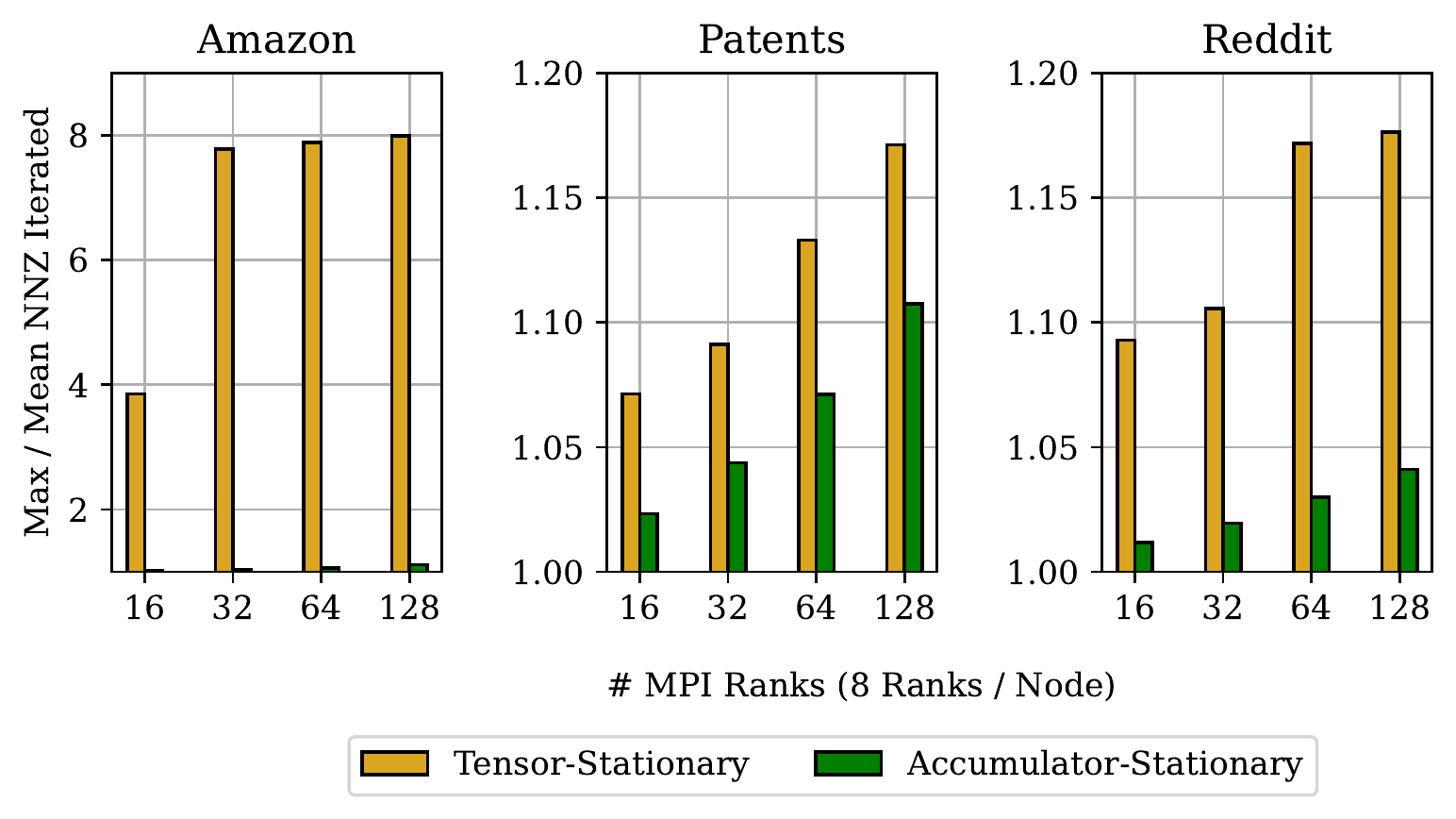}
    \caption{Avg. Load imbalance, defined
    as maximum / mean nonzeros iterated over by
    any MPI process, for d-STS-CP. 8 MPI 
    Ranks / Node, 5 Trials.}
    \Description[Load Imbalance Bar Charts]{Three bar 
    charts, one for Amazon, Patents, and Reddit, showing 
    load imbalance as the maximum / mean nonzeros iterated.
    Bars compare the tensor-stationary vs. 
    accumulator-stationary distribution at different MPI
    process counts. The tensor-stationary algorithm 
    has higher load imbalance across all tested 
    configurations.} 
    \label{fig:load_imbalance}
\end{figure}

Besides differences in the communication times
of the tensor-stationary and accumulator-stationary
schedules, Figure \ref{fig:communication_comparison} 
indicates a runtime difference in the downsampled MTTKRP
between the two schedules. 
Figure \ref{fig:load_imbalance} offers an 
explanation by comparing the load balance of these
methods. We measure load imbalance (averaged over 5 trials) as the maximum number of nonzeros processed in
the MTTKRP by any MPI process over the mean
of the same quantity. 

The accumulator-stationary schedule yields
better load balance over all tensors, with a dramatic
difference for the case of Amazon. The
latter exhibits a few rows of the Khatri-Rao
design matrix with high statistical leverage
and corresponding fibers with high nonzero
counts, producing the imbalance. The accumulator-stationary 
distribution (aided by the load balancing 
random permutation) distributes the nonzeros
in each selected fiber across all $P$ processors,
correcting the imbalance.

\subsection{Impact of Sample Count}

\begin{figure}
    \centering
    \includegraphics[scale=0.30]{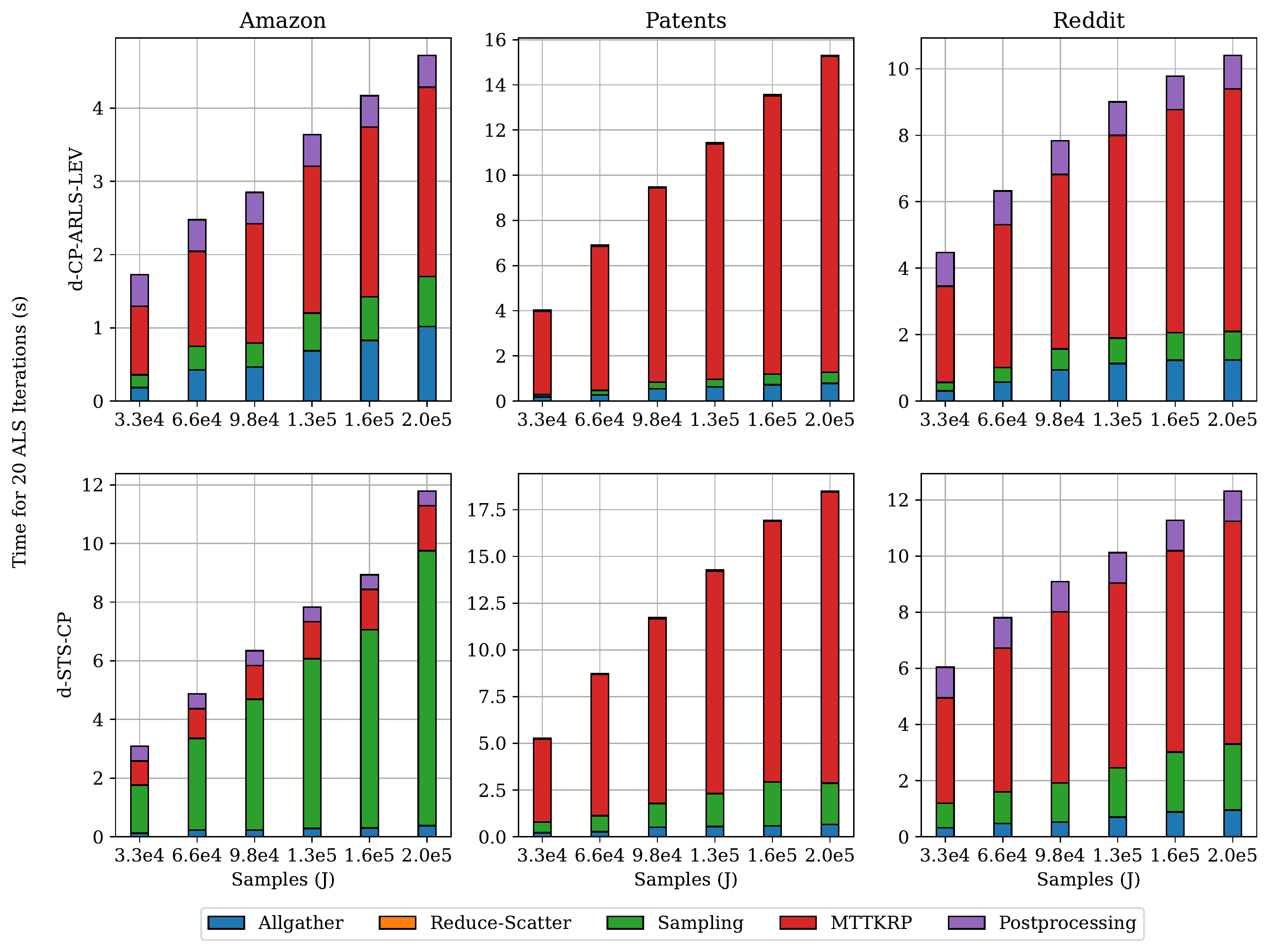}
    \caption{Runtime breakdown vs. sample count,
    $R=25$. 512 CPU cores, 5 trials, accumulator-stationary
    distribution.}
    \Description[Runtime breakdown for varying sample counts]
    {Two rows of three bar charts, one for Amazon, 
    Patents, and Reddit, showing the runtime breakdown 
    as the sample count increases from 32768 to 196608. 
    Top row shows d-CP-ARLS-LEV, bottom row shows d-STS-CP. 
    The runtime increase is roughly linear with 
    the sample count; for all charts but one, the 
    extra runtime is concentrated in the MTTKRP. For
    d-STS-CP on the Amazon tensor, runtime is concentrated in
    sample selection.} 
    \label{fig:sample_scaling}
\end{figure}

\begin{figure}
    \centering
    \includegraphics[scale=0.260]{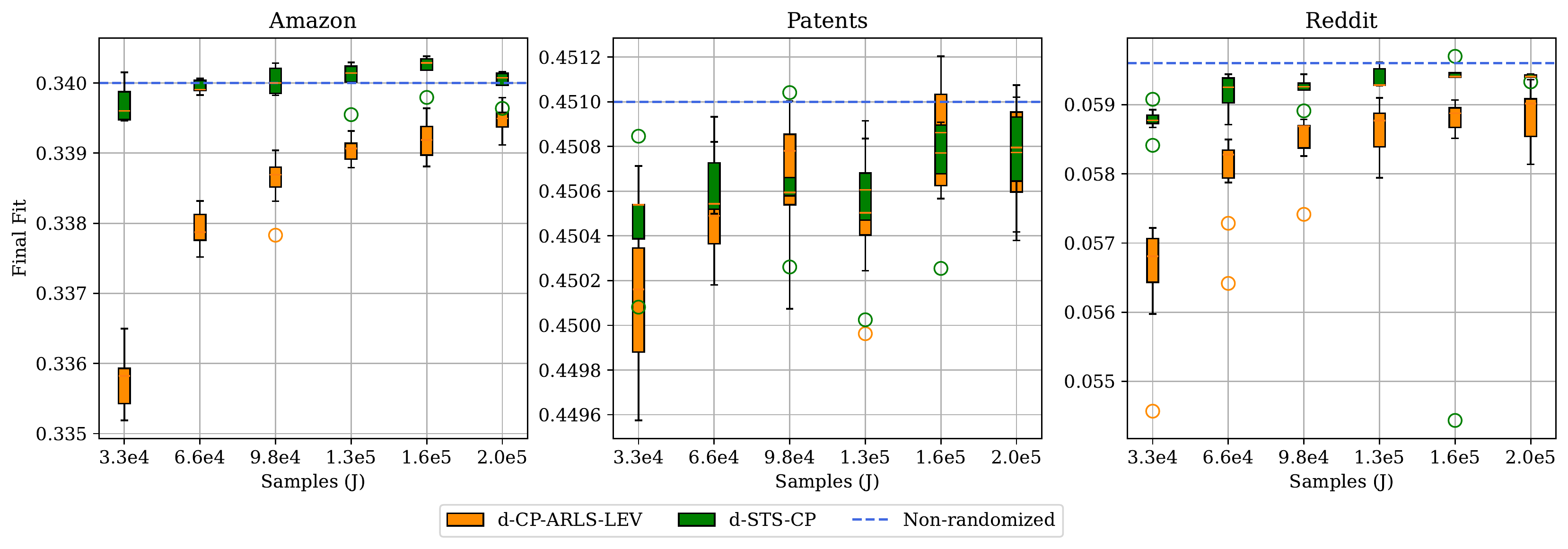}
    \caption{Final fit of randomized CP decomposition
    for varying sample count $J$. 
    Horizontal dashed lines indicate the
    fit produced by SPLATT. ALS was run for 
    40 iterations on Amazon and Patents, 
    80 iterations on Reddit, and for 10 
    trials each. All
    other experimental configuration 
    is identical to
    Figure \ref{fig:sample_scaling}.}
    \Description[Accuracy vs. Sample Count
    for Algorithms]
    {One row of three charts (one for Amazon,
    Patents, and Reddit). Each chart has
    boxplots showing the accuracy (y-axis) for
    varying sample counts. The accuracy
    for STS-CP is generally higher than
    CP-ARLS-LEV, but increasing the
    sample count does not produce significant
    accuracy gains at rank 25.} 
    \label{fig:sample_accuracy}
\end{figure}

In prior sections, we used the sample count $J=2^{16}$ to establish a consistent comparison
with prior work. Figure \ref{fig:sample_scaling} demonstrates
the runtime impact of increasing the sample count for both
of our algorithms on all three tensors. For all experiments
but one, the MTTKRP component of the runtime increases the
most as $J$ gets larger. For d-STS-CP on the Amazon tensor, the runtime increase owes primarily to the
higher cost of sample selection. The higher sampling
time for d-STS-CP on Amazon is explained in 
Section \ref{sec:strong_scaling}. 
Figure \ref{fig:sample_accuracy} gives the 
final fits after running our randomized 
algorithms for varying sample counts. The
increase in accuracy is minimal beyond
$J=2^{16}$ for d-STS-CP on Amazon and Reddit.
Both algorithms perform comparably on Patents.
These plots suggest that sample count as low
as $J=2^{16}$ is sufficient to achieve
competitive performance with libraries
like SPLATT on large tensors.

\section{Conclusions and Further Work}
We have demonstrated in this work that randomized 
CP decomposition algorithms are competitive at the scale
of thousands of CPU cores with state-of-the-art, 
highly-optimized non-randomized libraries for the same task. Future work includes improving the irregular 
communication pattern of the
d-STS-CP algorithm, as well as deploying our algorithm 
on massive real-world tensors larger than those
offered by FROSTT.

\section{Acknowledgements and Funding}
We thank the referees for valuable feedback which helped improve the paper. V. Bharadwaj was supported by the U.S. Department 
of Energy, Office of Science, Office of Advanced Scientific 
Computing Research, Department of Energy Computational Science 
Graduate Fellowship under Award Number DE-SC0022158. O. A. Malik and
A. Bulu\c{c} were supported by the Office of Science of the DOE
under Award Number DE-AC02-05CH11231. R. Murray was supported 
by Laboratory Directed Research and Development (LDRD) funding 
from Berkeley Lab, provided by the Director, Office of Science, of the U.S. DOE 
under Contract No. DE-AC02-05CH11231. R. Murray was also funded by an NSF Collaborative Research Framework under NSF Grant Nos. 2004235 and 2004763. This research used resources 
of the National Energy Research Scientific Computing Center, a DOE Office 
of Science User Facility, using 
NERSC award ASCR-ERCAP0024170. 

\appendix

\section{Distributed CP-ARLS-LEV Sampling}
Let $i_1, ..., i_{N-1}$ denote row indices from
factor matrices $U_1, ..., U_{N-1}$ 
that uniquely identify a row from the 
Khatri-Rao product $U_{\neq N}$. 
To efficiently sample according to
an \textit{approximate leverage score distribution}
on the rows of $U_{\neq N}$, the CP-ARLS-LEV 
algorithm by Larsen and Kolda \cite{larsen_practical_2022} weights each row by 
\[
\tilde \ell_{i_1, ..., i_{N-1}} := 
\prod_{k=1}^{N-1} U_k\br{i_k, :} 
G_k^+ U_k\br{i_k, :}^\top
\]
where $G_k := U_k^\top U_k$ for all $k$.
Because each weight in the distribution above is a
product of scores from each factor,
we can draw $i_1, ..., i_{N-1}$ independently 
and concatenate the indices to assemble
one row sample. Given that the
factors are distributed by 
block rows among processors, the main challenge is
to sample without gathering the
probability weight vector for each $U_k$ to
a single processor.

\begin{algorithm}
   \caption{CP-ARLS-LEV-build 
   $\paren{U_i^{(p_j)}}$)} 
   \label{alg:cp_arls_lev_build}
\begin{algorithmic}[1]
    \STATE $G_i := \textbf{Allreduce}
    \paren{U_i^{(p_j)\top} U_i^{(p_j)}}$

    \STATE $\textrm{dist}_i^{(p_j)} 
    := \textrm{diag}\paren{U_i^{(p_j)}G_i^+ U_i^{(p_j)\top}}$

    \STATE $C_i^{(p_j)} 
    := \norm{\textrm{dist}_i^{(p_j)}}_1$
    
    \STATE $\textrm{dist}_i^{(p_j)} /= C_i^{(p_j)}$

    \STATE \textbf{Postcondition: } 
    $G_i^+$, $\textrm{dist}_i^{(p_j)}$, and
    $C_i^{(p_j)}$ are initialized on each processor. 
\end{algorithmic}
\end{algorithm}

\label{sec:cp_arls_lev_sampling}
Algorithms \ref{alg:cp_arls_lev_build} and
\ref{alg:cp_arls_lev_sample} give full procedures
to build the distributed CP-ARLS-LEV 
data structure and draw samples from it, 
respectively. The
build algorithm is called for all $U_i$,
$1 \leq i \leq N$, before the ALS algorithm 
begins. It is also called each
time a matrix $U_i$ is updated in an ALS round.
Each procedure is executed
synchronously by all processors $p_j$, 
$1 \leq j \leq P$. Recall further that we 
define $U_i^{(p_j)}$ as the block row of 
the $i$-th factor matrix
uniquely owned by processor $p_j$. 
Algorithm \ref{alg:cp_arls_lev_build} allows 
all processors to redundantly compute the 
Gram matrix $G_i$ and the normalized 
local leverage score distribution on the 
block row $U_i^{(p_j)}$. 

Algorithm \ref{alg:cp_arls_lev_sample} enables 
each processor to
draw samples from the Khatri-Rao product 
$U_{\neq k}$. For each index
$i \neq k$, each processor determines 
the fraction of $J$ rows drawn from 
its local block using
a consistent multinomial sample
according to the weights $C_i^{(p_j)}$,
$1 \leq j \leq P$. By 
\textit{consistent}, we mean that
each processor executes the multinomial sampling
using a pseudorandom number generator with
a common seed that is shared
among all processors. The result of this
operation is a vector $SC^{\textrm{loc}} \in \ZZ^P$ 
which gives the sample count each processor should
draw locally. Each processor then samples 
from its local distribution. At the end of 
the algorithm, each
row of $X$ contains a sample drawn according
to the approximate leverage score distribution.

We note that the sampling algorithm, as 
presented, involves the \verb|Allgather| of
a $J \times N$ sampling matrix followed by
a random permutation. We use this procedure
in our code, since we found that the communication
cost $O(JN)$ was negligible for the range of
sample counts we used. However, this cost can
be reduced to $O(JN / P)$, in expectation,
with an \verb|All-to-allv| communication pattern
that permutes the indices without gathering them
to a single processor. 

\begin{algorithm}
   \caption{CP-ARLS-LEV-sample $\paren{k, J}$} 
   \label{alg:cp_arls_lev_sample}
\begin{algorithmic}[1]
    \STATE \textbf{Require: } Vectors 
    $\textrm{dist}_i^{(p_j)}$, and
    normalization constants 
    $C_i^{(p_j)}$.

    \STATE Initialize sample matrix 
    $X \in \ZZ^{J \times N}$ on all processors.

    \FOR{$i=1...N$, $i \neq k$}
        \STATE $C := \textbf{Allgather}\paren{C_i^{(p_j)}}$

        \STATE $W = \sum_{\ell=1}^P C\br{\ell}$

        \STATE $SC^{\textrm{loc}} :=
        \textrm{consistent-multinomial}(\br{C\br{1} / W, ..., C\br{P} / W}, J)$

        \STATE $\textrm{samples}^{\textrm{loc}}
        := 
        \textrm{sample}\paren{
        \textrm{dist}_i^{(p_j)}, SC^{\textrm{loc}}\br{j}}$

        \STATE $X\br{:, i} 
        := \textbf{Allgather}
        \paren{
        \textrm{samples}^{\textrm{loc}}
        }$\ \ \ \ //See note

        \STATE Perform a consistent
        random permutation of $X\br{:, i}$
    \ENDFOR
    \STATE \textbf{return} $X$, a set of
    samples from the Khatri-Rao product
    $U_{\neq k}$. 

\end{algorithmic}
\end{algorithm}

\section{Distributed STS-CP Sampling}
We use the same variables defined
at the beginning of Appendix 
\ref{sec:cp_arls_lev_sampling}. To draw samples
from the \textit{exact leverage score distribution},
the STS-CP algorithm conditions each row index draw
$i_k$ on draws $i_1, ..., i_{k-1}$
\cite{bharadwaj2023fast}. To formalize 
this, let $\hat i_1, ..., \hat i_{N-1}$ be random
variables for each index that jointly follow the
exact leverage distribution. Suppose we have
already sampled $\hat i_1 = i_1, ..., \hat i_{k-1} =
i_{k-1}$, and let 
$h = \startimes_{j=1}^{k-1} U_j\br{i_j, :}^\top$
be the product of these sampled rows. 
Bharadwaj et. al. show that the 
conditional probability of $\hat i_k = i_k$ is  
\[
p(\hat i_k = i_k\ \vert\ \hat i_{<k} = i_{<k}) \propto
\paren{U_k\br{i_k, :}^\top \circledast h}^\top G_{>k}
\paren{U_k\br{i_k, :}^\top \circledast h}
\]
where 
$G_{>k} = G^+ \circledast \startimes_{j=k}^{N-1} G_i$
\cite{bharadwaj2023fast}. The STS-CP algorithm exploits
this formula to efficiently sample from the exact
leverage distribution.

\label{sec:sts_cp_sampling}
Algorithms \ref{alg:sts_cp_build} and
\ref{alg:sts_cp_sample} give procedures to
build and sample from the distributed data 
structure for STS-CP, which are analogues of Algorithms
\ref{alg:cp_arls_lev_build} and 
\ref{alg:cp_arls_lev_sample} for CP-ARLS-LEV. To simplify 
our presentation, 
we assume that the processor count $P$ is 
a power of two. The general
case is a straightforward extension 
(see Chan et. al. \cite{chan_collective_2007}), and
our implementation makes no restriction on $P$. 
The build procedure in 
Algorithm \ref{alg:sts_cp_build} computes the
Gram matrix $G_i$ for each matrix $U_i$ using
a the bidirectional exchange algorithm for
Allreduce \cite{chan_collective_2007}. The 
difference is that each processor caches the
intermediate matrices that arise during the
reduction procedure, each uniquely identified 
with internal nodes of the binary tree in 
Figure \ref{fig:sts_cp_parallelization}.

\begin{algorithm}
   \caption{STS-CP-build 
   $\paren{U_i^{(p_j)}}$)} 
   \label{alg:sts_cp_build}
\begin{algorithmic}[1]
    \STATE $\tilde G_{\log_2 P} := 
    U_i^{(p_j)\top} U_i^{(p_j)}$

    \FOR{$\ell=\log_2 P...2$}
        \STATE \textbf{Send} $\tilde G_{\ell}$ to
        sibling of ancestor at level $\ell$, and \textbf{receive}
        the corresponding matrix 
        $\tilde G_{\textrm{sibling}}$.

        \STATE Assign $\tilde G_{\ell-1} =
        \tilde G_{\textrm{sibling}} + 
        \tilde G_{\ell}$

        \IF{Ancestor at level $\ell$ is a left child}
            \STATE $\tilde G_{\ell-1}^L := \tilde G_{\ell}$ 
        \ELSE 
            \STATE $\tilde G_{\ell-1}^L := \tilde G_{\textrm{sibling}}$ 
        \ENDIF
    \ENDFOR

    \STATE Assign $G_i := \tilde G_1$

    \STATE \textbf{Postcondition: } Each
    processor stores a list of \textit{partial
    gram matrices} $\tilde G_\ell$ and $G_\ell^L$, 
    from the root to its unique tree node. $G_i$
    is initialized.
\end{algorithmic}
\end{algorithm}

\begin{algorithm}
   \caption{STS-CP-sample 
   $\paren{k, J}$)} 
   \label{alg:sts_cp_sample}
\begin{algorithmic}[1]
    \STATE Initialize $X \in \ZZ^{J \times N}$, 
    $H \in \RR^{J \times (R+1)}$ distributed by block rows.
    Let $X^{(p_j)}, H^{(p_j)}$ be the block rows assigned to $p_j$.

    \STATE $X^{(p_j)} := \br{0}, H^{(p_j)} := \br{1}$
 
    \FOR{$i=1...N, i \neq k$}
        \STATE $G_{>k} := G^+ \circledast 
        \startimes_{\ell=k+1}^N G_i$ 
        \STATE $H^{(p_j)}\br{:, R+1} := 
        \textrm{uniform-samples}(\br{0, 1})$

        \FOR{$\ell=1...\log P - 1$}
            \STATE $J^{\textrm{loc}} = 
            \textrm{row-count}\paren{X^{(p_j)}}$
            \FOR{$k = 1...J^{\textrm{loc}}$}
                 \STATE $r := H^{(p_j)}\br{k, R+1}$ 
                 \STATE $h := H^{(p_j)}\br{k, 1:R}$
                 \STATE $X^{(p_j)}\br{i, k} *= 2$
                 \STATE $T = h^\top \paren{\tilde G_\ell^L \circledast G_{>k}} h / \paren{h^\top (\tilde 
                 G_\ell \circledast G_{>k}) h}$
                 \IF{$r \geq T$}
                    \STATE $X^{(p_j)}\br{i, k} \pluseq 1$
                    \STATE $r := (r - T) / (1 - T)$
                \ELSE
                    \STATE $r := r / T$
                \ENDIF 
                 \STATE $H^{(p_j)}\br{k, R+1} := r$ 
            \ENDFOR
            \STATE Execute an \textbf{All-to-allv} call to 
            redistribute $X^{(p_j)}$, $H^{(p_j)}$ according to the
            binary-tree data structure.
        \ENDFOR
        \STATE $J^{\textrm{loc}} = 
        \textrm{row-count}\paren{X^{(p_j)}}$
        \FOR{$k = 1...J^{\textrm{loc}}$}  
            \STATE $\textrm{idx} := 
            \textrm{local-STS-CP}(H^{(p_j)}\br{k, 1:R}, \tilde 
            G_{\log 2 P}, G_{>k}, r)$  
            \STATE $X^{(p_j)}\br{i, k} := \textrm{idx}$
            
            \STATE $H^{(p_j)}\br{k, 1:R} *= U_i^{(p_j)}\br{\textrm{idx} - I_i p_j / P, :}$ 
        \ENDFOR
    \ENDFOR
    \STATE \textbf{return} $X^{(p_j)}, H^{(p_j)}$
\end{algorithmic}
\end{algorithm}

In the sampling algorithm, the cached matrices are used to 
determine the index of a row drawn from $U_i$ via binary search.
The matrix of sample indices $X$ and sampled rows $H$ are initially
distributed by block rows among processors. Then for each matrix
$U_i$, $i \neq k$, a random number is drawn uniformly in the
interval $\br{0,1}$ for each sample. By stepping down levels of
the tree, $J$ binary searches are computed in parallel to determine
the containing bin of each random draw. At each level, the cached
matrices tell the program whether to branch left or right
by computing the branching threshold $T$, which is compared
to the random draw $r$. The values in each column of 
$X^{(p_j)}$ hold the current node index of each sample 
at level $\ell$ of the search.
At level $L=\log_2 P$, the algorithm continues 
the binary search locally on each processor until 
a row index is identified (a procedure we denote as
``local-STS-CP". For more details, see 
the original work \cite{bharadwaj2023fast}. At 
the end of the algorithm, the sample indices 
in $X$ are correctly drawn according
to the exact leverage scores of $U_\neq k$.
The major communication cost of this algorithm stems from the
\textbf{All-to-allv} collective between levels of the
binary search. Because a processor may not have the 
required matrices $\tilde G_\ell, \tilde G^L_\ell$ to
compute the branching threshold for a sample, the sample must
be routed to another processor that owns the information. 
\newpage

\bibliographystyle{IEEEtran}
\balance
\bibliography{main}

\begin{thebibliography}{10}
\providecommand{\url}[1]{#1}
\csname url@samestyle\endcsname
\providecommand{\newblock}{\relax}
\providecommand{\bibinfo}[2]{#2}
\providecommand{\BIBentrySTDinterwordspacing}{\spaceskip=0pt\relax}
\providecommand{\BIBentryALTinterwordstretchfactor}{4}
\providecommand{\BIBentryALTinterwordspacing}{\spaceskip=\fontdimen2\font plus
\BIBentryALTinterwordstretchfactor\fontdimen3\font minus
  \fontdimen4\font\relax}
\providecommand{\BIBforeignlanguage}[2]{{%
\expandafter\ifx\csname l@#1\endcsname\relax
\typeout{** WARNING: IEEEtran.bst: No hyphenation pattern has been}%
\typeout{** loaded for the language `#1'. Using the pattern for}%
\typeout{** the default language instead.}%
\else
\language=\csname l@#1\endcsname
\fi
#2}}
\providecommand{\BIBdecl}{\relax}
\BIBdecl

\bibitem{smith_frostt_2017}
\BIBentryALTinterwordspacing
S.~Smith, J.~W. Choi, J.~Li, R.~Vuduc, J.~Park, X.~Liu, and G.~Karypis,
  ``{FROSTT}: {The} {Formidable} {Repository} of {Open} {Sparse} {Tensors} and
  {Tools},'' 2017. [Online]. Available: \url{http://frostt.io/}
\BIBentrySTDinterwordspacing

\bibitem{mao_malspot_2014}
H.-H. Mao, C.-J. Wu, E.~E. Papalexakis, C.~Faloutsos, K.-C. Lee, and T.-C. Kao,
  ``\BIBforeignlanguage{en}{{MalSpot}: {Multi2} {Malicious} {Network}
  {Behavior} {Patterns} {Analysis}},'' in
  \emph{\BIBforeignlanguage{en}{Advances in {Knowledge} {Discovery} and {Data}
  {Mining}}}, ser. Lecture {Notes} in {Computer} {Science}, V.~S. Tseng, T.~B.
  Ho, Z.-H. Zhou, A.~L.~P. Chen, and H.-Y. Kao, Eds.\hskip 1em plus 0.5em minus
  0.4em\relax Cham: Springer International Publishing, 2014, pp. 1--14.

\bibitem{balazevic_tucker_2019}
I.~Balazevic, C.~Allen, and T.~Hospedales, ``{TuckER}: {Tensor} {Factorization}
  for {Knowledge} {Graph} {Completion},'' in \emph{Proceedings of the 2019
  {Conference} on {Empirical} {Methods} in {Natural} {Language} {Processing}
  and the 9th {International} {Joint} {Conference} on {Natural} {Language}
  {Processing} ({EMNLP}-{IJCNLP})}.\hskip 1em plus 0.5em minus 0.4em\relax Hong
  Kong, China: Association for Computational Linguistics, Nov. 2019, pp.
  5185--5194.

\bibitem{haesun_matrix_nmf}
H.~Kim and H.~Park, ``{Sparse non-negative matrix factorizations via
  alternating non-negativity-constrained least squares for microarray data
  analysis},'' \emph{Bioinformatics}, vol.~23, no.~12, pp. 1495--1502, 05 2007.

\bibitem{gcp_kolda}
D.~Hong, T.~G. Kolda, and J.~A. Duersch, ``Generalized canonical polyadic
  tensor decomposition,'' \emph{SIAM Review}, vol.~62, no.~1, pp. 133--163,
  2020.

\bibitem{larsen_practical_2022}
B.~W. Larsen and T.~G. Kolda, ``Practical leverage-based sampling for low-rank
  tensor decomposition,'' \emph{SIAM J. Matrix Analysis and Applications}, June
  2022, accepted for publication.

\bibitem{smith_streaming}
S.~Smith, K.~Huang, N.~D. Sidiropoulos, and G.~Karypis, \emph{Streaming Tensor
  Factorization for Infinite Data Sources}.\hskip 1em plus 0.5em minus
  0.4em\relax {SIAM}, 2018, pp. 81--89.

\bibitem{kolda_tensor_2009}
T.~G. Kolda and B.~W. Bader, ``Tensor {Decompositions} and {Applications},''
  \emph{SIAM Review}, vol.~51, no.~3, pp. 455--500, Aug. 2009, publisher:
  Society for Industrial and Applied Mathematics.

\bibitem{choi_dfacto_2014}
J.~H. Choi and S.~Vishwanathan, ``{DFacTo}: {Distributed} {Factorization} of
  {Tensors},'' in \emph{Advances in {Neural} {Information} {Processing}
  {Systems}}, Z.~Ghahramani, M.~Welling, C.~Cortes, N.~Lawrence, and K.~Q.
  Weinberger, Eds., vol.~27.\hskip 1em plus 0.5em minus 0.4em\relax Curran
  Associates, Inc., 2014.

\bibitem{smith_splatt_2015}
S.~Smith, N.~Ravindran, N.~D. Sidiropoulos, and G.~Karypis, ``{SPLATT}:
  {Efficient} and {Parallel} {Sparse} {Tensor}-{Matrix} {Multiplication},'' in
  \emph{2015 {IEEE} {International} {Parallel} and {Distributed} {Processing}
  {Symposium}}, May 2015, pp. 61--70, iSSN: 1530-2075.

\bibitem{hypertensor}
O.~Kaya and B.~Uçar, ``Scalable sparse tensor decompositions in distributed
  memory systems,'' in \emph{SC '15: Proceedings of the International
  Conference for High Performance Computing, Networking, Storage and Analysis},
  2015, pp. 1--11.

\bibitem{bigtensor_16}
N.~Park, B.~Jeon, J.~Lee, and U.~Kang, ``Bigtensor: Mining billion-scale tensor
  made easy,'' in \emph{Proceedings of the 25th ACM International on Conference
  on Information and Knowledge Management}, ser. CIKM '16.\hskip 1em plus 0.5em
  minus 0.4em\relax New York, NY, USA: Association for Computing Machinery,
  2016, p. 2457–2460.

\bibitem{cheng_spals_2016}
D.~Cheng, R.~Peng, Y.~Liu, and I.~Perros, ``{SPALS}: {Fast} {Alternating}
  {Least} {Squares} via {Implicit} {Leverage} {Scores} {Sampling},'' in
  \emph{Advances in {Neural} {Information} {Processing} {Systems}}, D.~Lee,
  M.~Sugiyama, U.~Luxburg, I.~Guyon, and R.~Garnett, Eds., vol.~29.\hskip 1em
  plus 0.5em minus 0.4em\relax Curran Associates, Inc., 2016.

\bibitem{malik_more_2022}
O.~A. Malik, ``\BIBforeignlanguage{en}{More {Efficient} {Sampling} for {Tensor}
  {Decomposition} {With} {Worst}-{Case} {Guarantees}},'' in
  \emph{\BIBforeignlanguage{en}{Proceedings of the 39th {International}
  {Conference} on {Machine} {Learning}}}.\hskip 1em plus 0.5em minus
  0.4em\relax PMLR, Jun. 2022, pp. 14\,887--14\,917, iSSN: 2640-3498.

\bibitem{bharadwaj2023fast}
\BIBentryALTinterwordspacing
V.~Bharadwaj, O.~A. Malik, R.~Murray, L.~Grigori, A.~Buluc, and J.~Demmel,
  ``{Fast} {Exact} {Leverage} {Score} {Sampling} from {Khatri}-{Rao} {Products}
  with {Applications} to {Tensor} {Decomposition},'' in \emph{Thirty-seventh
  Conference on Neural Information Processing Systems}, 2023. [Online].
  Available: \url{https://arxiv.org/pdf/2301.12584.pdf}
\BIBentrySTDinterwordspacing

\bibitem{smith_medium-grained_2016}
S.~Smith and G.~Karypis, ``A {Medium}-{Grained} {Algorithm} for {Sparse}
  {Tensor} {Factorization},'' in \emph{2016 {IEEE} {International} {Parallel}
  and {Distributed} {Processing} {Symposium} ({IPDPS})}, May 2016, pp.
  902--911, iSSN: 1530-2075.

\bibitem{mahoney_rand_survey}
M.~W. Mahoney, ``Randomized algorithms for matrices and data,''
  \emph{Foundations and Trends® in Machine Learning}, vol.~3, no.~2, pp.
  123--224, 2011.

\bibitem{drineas_mahoney_rnla_survey}
P.~Drineas and M.~W. Mahoney, ``{RandNLA}: Randomized numerical linear
  algebra,'' \emph{Commun. ACM}, vol.~59, no.~6, p. 80–90, may 2016.

\bibitem{martinsson_tropp_2020}
P.-G. Martinsson and J.~A. Tropp, ``Randomized numerical linear algebra:
  Foundations and algorithms,'' \emph{Acta Numerica}, vol.~29, p. 403–572,
  2020.

\bibitem{nisa_csf_mm}
I.~Nisa, J.~Li, A.~Sukumaran-Rajam, P.~S. Rawat, S.~Krishnamoorthy, and
  P.~Sadayappan, ``An efficient mixed-mode representation of sparse tensors,''
  in \emph{Proceedings of the International Conference for High Performance
  Computing, Networking, Storage and Analysis}, ser. SC '19.\hskip 1em plus
  0.5em minus 0.4em\relax New York, NY, USA: Association for Computing
  Machinery, 2019.

\bibitem{ahle_treesketch}
T.~D. Ahle, M.~Kapralov, J.~B.~T. Knudsen, R.~Pagh, A.~Velingker, D.~P.
  Woodruff, and A.~Zandieh, ``Oblivious sketching of high-degree polynomial
  kernels,'' in \emph{Proceedings of the Thirty-First Annual ACM-SIAM Symposium
  on Discrete Algorithms}, ser. SODA '20.\hskip 1em plus 0.5em minus
  0.4em\relax USA: Society for Industrial and Applied Mathematics, 2020, p.
  141–160.

\bibitem{practical_randomized_cp}
C.~Battaglino, G.~Ballard, and T.~G. Kolda, ``A practical randomized cp tensor
  decomposition,'' \emph{SIAM Journal on Matrix Analysis and Applications},
  vol.~39, no.~2, pp. 876--901, 2018.

\bibitem{gittens_adaptive_sketching}
A.~Gittens, K.~Aggour, and B.~Yener, ``Adaptive sketching for fast and
  convergent canonical polyadic decomposition,'' in \emph{Proceedings of the
  37th International Conference on Machine Learning}, ser. Proceedings of
  Machine Learning Research, vol. 119.\hskip 1em plus 0.5em minus 0.4em\relax
  PMLR, 13--18 Jul 2020, pp. 3566--3575.

\bibitem{smith_csf}
S.~Smith and G.~Karypis, ``Tensor-matrix products with a compressed sparse
  tensor,'' in \emph{Proceedings of the 5th Workshop on Irregular Applications:
  Architectures and Algorithms}, ser. IA<sup>3</sup> '15.\hskip 1em plus 0.5em
  minus 0.4em\relax New York, NY, USA: Association for Computing Machinery,
  2015.

\bibitem{parti}
\BIBentryALTinterwordspacing
J.~Li, Y.~Ma, and R.~Vuduc, ``{ParTI!} : A parallel tensor infrastructure for
  multicore cpus and gpus,'' Oct 2018, last updated: Jan 2020. [Online].
  Available: \url{http://parti-project.org}
\BIBentrySTDinterwordspacing

\bibitem{oom_sparse_mttkrp}
A.~Nguyen, A.~E. Helal, F.~Checconi, J.~Laukemann, J.~J. Tithi, Y.~Soh,
  T.~Ranadive, F.~Petrini, and J.~W. Choi, ``Efficient, out-of-memory sparse
  mttkrp on massively parallel architectures,'' in \emph{Proceedings of the
  36th ACM International Conference on Supercomputing}, ser. ICS '22.\hskip 1em
  plus 0.5em minus 0.4em\relax New York, NY, USA: Association for Computing
  Machinery, 2022.

\bibitem{phipps2019software}
E.~T. Phipps and T.~G. Kolda, ``Software for sparse tensor decomposition on
  emerging computing architectures,'' \emph{SIAM Journal on Scientific
  Computing}, vol.~41, no.~3, pp. C269--C290, 2019.

\bibitem{wijeratne2023dynasor}
S.~Wijeratne, R.~Kannan, and V.~Prasanna, ``Dynasor: A dynamic memory layout
  for accelerating sparse mttkrp for tensor decomposition on multi-core cpu,''
  2023.

\bibitem{kanakagiri2023minimum}
R.~Kanakagiri and E.~Solomonik, ``Minimum cost loop nests for contraction of a
  sparse tensor with a tensor network,'' 2023.

\bibitem{ballard_parallel_2018}
G.~Ballard, K.~Hayashi, and K.~Ramakrishnan, ``Parallel nonnegative {CP}
  decomposition of dense tensors,'' in \emph{2018 {IEEE} 25th {International}
  {Conference} on {High} {Performance} {Computing} ({HiPC})}.\hskip 1em plus
  0.5em minus 0.4em\relax IEEE, 2018, pp. 22--31.

\bibitem{kang_gigatensor_2012}
U.~Kang, E.~Papalexakis, A.~Harpale, and C.~Faloutsos, ``{GigaTensor}: scaling
  tensor analysis up by 100 times - algorithms and discoveries,'' in
  \emph{Proceedings of the 18th {ACM} {SIGKDD} international conference on
  {Knowledge} discovery and data mining}, ser. {KDD} '12.\hskip 1em plus 0.5em
  minus 0.4em\relax New York, NY, USA: Association for Computing Machinery,
  Aug. 2012, pp. 316--324.

\bibitem{ma_efficient_2021}
L.~Ma and E.~Solomonik, ``Efficient parallel {CP} decomposition with pairwise
  perturbation and multi-sweep dimension tree,'' in \emph{2021 {IEEE}
  {International} {Parallel} and {Distributed} {Processing} {Symposium}
  ({IPDPS})}, May 2021, pp. 412--421, iSSN: 1530-2075.

\bibitem{solomonik_massively_2014}
E.~Solomonik, D.~Matthews, J.~R. Hammond, J.~F. Stanton, and J.~Demmel, ``A
  massively parallel tensor contraction framework for coupled-cluster
  computations,'' \emph{Journal of Parallel and Distributed Computing},
  vol.~74, no.~12, pp. 3176--3190, 2014, publisher: Academic Press.

\bibitem{yadav_spdistal_2022}
R.~Yadav, A.~Aiken, and F.~Kjolstad, ``Spdistal: Compiling distributed sparse
  tensor computations,'' in \emph{Proceedings of the International Conference
  on High Performance Computing, Networking, Storage and Analysis}, ser. SC
  '22.\hskip 1em plus 0.5em minus 0.4em\relax IEEE Press, 2022.

\bibitem{jin_fjlt}
R.~Jin, T.~G. Kolda, and R.~Ward, ``{Faster Johnson–Lindenstrauss transforms
  via Kronecker products},'' \emph{Information and Inference: A Journal of the
  IMA}, vol.~10, no.~4, pp. 1533--1562, 10 2020.

\bibitem{diao_kronecker_sketch}
H.~Diao, Z.~Song, W.~Sun, and D.~Woodruff, ``Sketching for kronecker product
  regression and p-splines,'' in \emph{International Conference on Artificial
  Intelligence and Statistics}.\hskip 1em plus 0.5em minus 0.4em\relax PMLR,
  2018, pp. 1299--1308.

\bibitem{kolda_stochastic_2020}
T.~G. Kolda and D.~Hong, ``\BIBforeignlanguage{en}{Stochastic {Gradients} for
  {Large}-{Scale} {Tensor} {Decomposition}},''
  \emph{\BIBforeignlanguage{en}{SIAM Journal on Mathematics of Data Science}},
  vol.~2, no.~4, pp. 1066--1095, Jan. 2020.

\bibitem{chan_collective_2007}
E.~Chan, M.~Heimlich, A.~Purkayastha, and R.~van~de Geijn, ``Collective
  communication: theory, practice, and experience: {Research} {Articles},''
  \emph{Concurrency and Computation: Practice \& Experience}, vol.~19, no.~13,
  pp. 1749--1783, Sep. 2007.

\bibitem{rolinger_performance}
T.~B. Rolinger, T.~A. Simon, and C.~D. Krieger, ``Performance considerations
  for scalable parallel tensor decomposition,'' \emph{Journal of Parallel and
  Distributed Computing}, vol. 129, pp. 83--98, 2019.

\bibitem{aktulga_tall_skinny}
H.~M. Aktulga, A.~Buluç, S.~Williams, and C.~Yang, ``Optimizing sparse
  matrix-multiple vectors multiplication for nuclear configuration interaction
  calculations,'' in \emph{2014 IEEE 28th International Parallel and
  Distributed Processing Symposium}, 2014, pp. 1213--1222.

\end{thebibliography}

\end{document}